\documentclass[10pt]{amsart}
\usepackage{amsmath,amscd,amssymb,amsthm,amsbsy,amscd,stmaryrd}
\usepackage{mathrsfs}
\input{xy}

\usepackage{pstricks}
\usepackage{color}
\usepackage{graphicx}
\usepackage[all]{xy}
\usepackage {hyperref}

\newtheorem{theorem}{Theorem}[section]
\newtheorem{lemma}[theorem]{Lemma}
\newtheorem{proposition}[theorem]{Proposition}
\newtheorem{corollary}[theorem]{Corollary}
\newtheorem{remark}[theorem]{Remark}

\theoremstyle{definition}
\newtheorem{definition}[theorem]{Definition}

\numberwithin{equation}{section}

\begin{document}

\title{Eigenvalues and entropies under the harmonic-Ricci flow}
\author{Yi Li}
\address{Department of Mathematics, Harvard University, One Oxford street, Cambridge, MA 02138}

\curraddr{Department of Mathematics, Johns Hopkins University, 3400 N. Charles Street, Baltimore, MD, 21218}

\email{yli@math.jhu.edu}




\begin{abstract} In this paper, the author discusses the eigenvalues and entropies under the harmonic-Ricci flow, which is the Ricci flow coupled with the harmonic map flow. We give an alternative proof of results for compact steady and expanding harmonic-Ricci breathers. In the second part, we derive some monotonicity formulas for eigenvalues of Laplacian under the harmonic-Ricci flow. Finally, we obtain the first variation of the shrinker and expanding entropies of the harmonic-Ricci flow.
\end{abstract}

\maketitle

\tableofcontents

\section{Introduction}\label{section1}

After successfully applying the Ricci flow to topological and geometric problems, people study some analogue flows, including the harmonic-Ricci flow \cite{L,M2}, connection Ricci flow \cite{S2}, Ricci-Yang-Mills flow \cite{S1,S4,Y}, and renormalization group flows \cite{HHKL,LY,OSW,S3}, etc. In this note, we study the eigenvalue problems of the harmonic-Ricci flow which is the following coupled system
\begin{eqnarray}
\frac{\partial}{\partial t}g(t)&=&-2{\rm Ric}_{g(t)}+4du(t)\otimes du(t),\label{1.1}\\
\frac{\partial}{\partial t}u(t)&=&\Delta_{g(t)}u(t).\label{1.2}
\end{eqnarray}
For convenience, we introduce a new symmetric $2$-tensor $\mathcal{S}_{g(t),u(t)}$ whose components $S_{ij}$ are defined by
\begin{equation*}
S_{ij}:=R_{ij}-2\partial_{i}u\partial_{j}u.
\end{equation*}
Its trace is $S_{g(t),u(t)}:=g^{ij}S_{ij}=R_{g(t)}-2\left|\nabla_{g(t)} u(t)\right|^{2}_{g(t)}$.

Suppose that $M$ is a compact Riemannian manifold. For any
Riemannian metric $g$ and any smooth functions $u,f$, we have a number of functionals
\begin{eqnarray*}
\mathcal{F}(g,u,f)&=&\int_{M}\left(R_{g}+\left|\nabla_{g} f\right|^{2}_{g}
-2\left|\nabla_{g} u\right|^{2}_{g}\right)e^{-f}\!\ dV_{g},\\
\mathcal{E}(g,u,f)&=&\int_{M}\left(R_{g}-2\left|\nabla_{g} u\right|^{2}_{g}\right)e^{-f}\!\ dV_{g},\\
\mathcal{F}_{k}(g,u,f)&=&\int_{M}\left(kR_{g}+\left|\nabla_{g} f\right|^{2}_{g}
-2k\left|\nabla_{g} u\right|^{2}_{g}\right)e^{-f}\!\ dV_{g}.
\end{eqnarray*}
List \cite{L} and M\"uller \cite{M2} showed that, as in the case of Perelman's $\mathcal{F}$-functional, under the following evolution equation
\begin{eqnarray}
\frac{\partial}{\partial t}g(t)&=&-2{\rm Ric}_{g(t)}+4du(t)\otimes du(t),\nonumber\\
\frac{\partial}{\partial t}u(t)&=&\Delta_{g(t)}u(t),\label{1.3}\\
\frac{\partial}{\partial t}f(t)&=&-\Delta_{g(t)}f(t)-R_{g(t)}+\left|\nabla_{g(t)} f(t)\right|^{2}_{g(t)}+2\left|\nabla_{g(t)} u(t)\right|^{2}_{g(t)},\nonumber
\end{eqnarray}
the evolution equation for $\mathcal{F}$-functional is
\begin{eqnarray}
\frac{d}{dt}\mathcal{F}(g(t),u(t),f(t))&=&2\int_{M}
\left|\mathcal{S}_{g(t),u(t)}+\nabla^{2}_{g(t)}f(t)\right|^{2}_{g(t)}
e^{-f(t)}\!\ dV_{g(t)}\nonumber\\
&&+ \ 4\int_{M}\left|\Delta_{g(t)}u(t)-\langle du(t),df(t)\rangle_{g(t)}\right|^{2}_{g(t)}e^{-f(t)}\!\ dV_{g(t)}\label{1.4}
\end{eqnarray}
that is nonnegative. Based on (\ref{1.4}), we derive

\begin{theorem}\label{t1.1} Under the evolution equation (\ref{1.3}), one has
\begin{eqnarray}
\frac{d}{dt}\mathcal{E}(g(t),u(t),f(t))&=&2\int_{M}\left|\mathcal{S}_{g(t),u(t)}\right|^{2}_{g(t)}
e^{-f(t)}\!\ dV_{g(t)}\nonumber\\
&&+ \ 4\int_{M}\left|\Delta_{g(t)}u(t)\right|^{2}_{g(t)}e^{-f(t)}\!\ dV_{g(t)},\label{1.5}\\
\frac{d}{dt}\mathcal{F}_{k}(g(t),u(t),f(t))&=&2(k-1)\int_{M}
\left|\mathcal{S}_{g(t),u(t)}\right|^{2}_{g(t)}
e^{-f(t)}\!\ dV_{g(t)}\nonumber\\
&&+2\int_{M}\left|\mathcal{S}_{g(t),u(t)}+\nabla^{2}_{g(t)}
f(t)\right|^{2}_{g(t)}e^{-f(t)}\!\ dV_{g(t)}\label{1.6}\\
&&+ \ 4(k-1)\int_{M}\left|\Delta_{g(t)}u(t)\right|^{2}_{g(t)}e^{-f(t)}\!\ dV_{g(t)}\nonumber\\
&&+ \ 4\int_{M}\left|\Delta_{g(t)}u(t)-\langle du(t),df(t)\rangle_{g(t)}\right|^{2}_{g(t)}e^{-f(t)}\!\ dV_{g(t)}\nonumber.
\end{eqnarray}
\end{theorem}

As a corollary we give a new proof of the following

\begin{corollary}\label{c1.2} There is no compact steady harmonic-Ricci breather other than $(M,g(t))$ is Ricci-flat and $u(t)$ is constant.
\end{corollary}

When we deal with the expanding harmonic-Ricci breather, we need the following two functionals
\begin{eqnarray*}
\mathcal{L}_{+}(g,u,\tau,f)&=&\tau^{2}\int_{M}\left(R_{g}+\frac{n}{2\tau}
+\Delta_{g}f-2\left|\nabla_{g} u\right|^{2}_{g}\right)e^{-f}\!\ dV_{g},\\
\mathcal{L}_{+,k}(g,u,\tau,f)&=&\tau^{2}\int_{M}\left[k
\left(R_{g}+\frac{n}{2\tau}\right)
+\Delta_{g}f-2k\left|\nabla_{g} u\right|^{2}_{g}\right]e^{-f}\!\ dV_{g}.
\end{eqnarray*}

Under the following evolution equation
\begin{eqnarray*}
\frac{\partial}{\partial t}g(t)&=&-2{\rm Ric}_{g(t)}+4du(t)\otimes du(t),\\
\frac{\partial}{\partial t}u(t)&=&\Delta_{g(t)}u(t),\\
\frac{\partial}{\partial t}f(t)&=&-\Delta_{g(t)}f(t)+\left|\nabla_{g(t)} f(t)\right|^{2}_{g(t)}-R_{g(t)}+2\left|\nabla_{g(t)} u(t)\right|^{2}_{g(t)},\\
\frac{d}{dt}\tau(t)&=&1,
\end{eqnarray*}
we have

\begin{theorem}\label{t1.3} Under the above evolution equation, one has
\begin{eqnarray}
&&\frac{d}{dt}\mathcal{L}_{+}(g(t),u(t),\tau(t),f(t))\nonumber\\
&=&2\tau(t)^{2}\int_{M}\left|\mathcal{S}_{g(t),u(t)}
+\nabla^{2}_{g(t)}f(t)+\frac{1}{2\tau(t)}g(t)\right|^{2}_{g(t)}
e^{-f(t)}\!\ dV_{g(t)}\label{1.7}\\
&&+4\tau(t)^{2}\int_{M}\left|\Delta_{g(t)}u(t)-\langle du(t),df(t)\rangle_{g(t)}\right|^{2}_{g(t)}e^{-f(t)}\!\ dV_{g(t)},\nonumber\\
&&\frac{d}{dt}\mathcal{L}_{+,k}(g(t),u(t),\tau(t),f(t))\nonumber\\
&=&2\tau(t)^{2}\int_{M}\left|\mathcal{S}_{g(t),u(t)}
+\nabla^{2}_{g(t)}f(t)+\frac{1}{2\tau(t)}g(t)\right|^{2}_{g(t)}
e^{-f(t)}dV_{g(t)}\nonumber\\
&&+ \ 2(k-1)\tau(t)^{2}\int_{M}\left|\mathcal{S}_{g(t),u(t)}+\frac{1}{2\tau(t)}g(t)\right|^{2}_{g(t)}
e^{-f(t)}dV_{g(t)}\label{1.8}\\
&&+ \ 4\tau(t)^{2}\int_{M}\left|\Delta_{g(t)}u(t)-\langle du(t),df(t)\rangle_{g(t)}\right|^{2}_{g(t)}e^{-f(t)}dV_{g(t)}\nonumber\\
&&+ \ 4(k-1)\tau(t)^{2}\int_{M}\left|\Delta_{g(t)}u(t)\right|^{2}_{g(t)}
e^{-f(t)}dV_{g(t)}.\nonumber
\end{eqnarray}
\end{theorem}

As a corollary, we obtain a new proof of the following

\begin{corollary} \label{c1.4}There is no expanding harmonic-Ricci breather on compact Riemannian manifolds other than $M$ is an Einstein manifold and $u(t)$ is constant.
\end{corollary}

The second part of this paper focuses on the eigenvalue of the Laplacian operator under the harmonic-Ricci flow.

\begin{theorem}\label{t1.5} If $(g(t),u(t))$ is a solution of the harmonic-Ricci flow on a compact Riemannian manifold $M$ and $\lambda(t)$ denotes the eigenvalue of the Laplacian $\Delta_{g(t)}$ with eigenfunction $f(t)$, then
\begin{eqnarray}
\frac{d}{dt}\lambda(t)\cdot\int_{M}f(t)^{2}\!\ dV_{g(t)}&=&\lambda(t)\int_{M}S_{g(t),u(t)}f(t)^{2}\!\ dV_{g(t)}\nonumber\\
&&- \ \int_{M}S_{g(t),u(t)}\left|\nabla_{g(t)} f\right|^{2}_{g(t)}\!\ dV_{g(t)}\label{1.9}\\
&&+ \ 2\int_{M}\langle\mathcal{S}_{g(t),u(t)},df(t)\otimes df(t)\rangle_{g(t)}\!\ dV_{g(t)}.\nonumber
\end{eqnarray}
\end{theorem}

The above equation (\ref{1.9}) is a general formula to describe the evolution of $\lambda(t)$ under the harmonic-Ricci flow. Under a curvature assumption, we can derive some monotonicity formulas for the eigenvalue $\lambda(t)$. Set
\begin{equation}
S_{\min}(0):=\min_{x\in M}S_{g(0),u(0)}(x)\label{1.10}
\end{equation}
the minimum of $S_{g(t),u(t)}$ over $M$ at the time $0$.

\begin{theorem}\label{t1.6} Let $(g(t),u(t))_{t\in[0,T]}$ be a solution of the harmonic-Ricci flow on a compact Riemannian manifold $M$ and $\lambda(t)$ denote the eigenvalue of the Laplacian $\Delta_{g(t)}$. Suppose that $\mathcal{S}_{g(t),u(t)}-\alpha S_{g(t),u(t)}g(t)\geq0$ along the harmonic-Ricci flow for some $\alpha\geq\frac{1}{2}$.

\begin{itemize}

\item[(1)] If $S_{\min}(0)\geq0$, then $\lambda(t)$ is nondecreasing along the harmonic-Ricci flow for any $t\in[0,T]$.

\item[(2)] If $S_{\min}(0)>0$, then the quantity
\begin{equation*}
\left(1-\frac{2}{n}S_{\min}(0)t\right)^{n\alpha}\lambda(t)
\end{equation*}
is nondecreasing along the harmonic-Ricci flow for
$T\leq\frac{n}{2S_{\min}(0)}$.

\item[(3)] If $S_{\min}(0)<0$, then the quantity
\begin{equation*}
\left(1-\frac{2}{n}S_{\min}(0)t\right)^{n\alpha}\lambda(t)
\end{equation*}
is nondecreasing along the harmonic-Ricci flow for any $t\in[0,T]$.

\end{itemize}
\end{theorem}

\begin{corollary}\label{c1.7} Let $(g(t),u(t))_{t\in[0,T]}$ be a solution of the harmonic-Ricci flow on a compact Riemannian surface $\Sigma$ and $\lambda(t)$ denote the eigenvalue of the Laplacian $\Delta_{g(t)}$.

\begin{itemize}

\item[(1)] Suppose that ${\rm Ric}_{g(t)}\leq\epsilon du(t)\otimes du(t)$ where
\begin{equation*}
\epsilon\leq4\frac{1-\alpha}{1-2\alpha}, \ \ \ \alpha>\frac{1}{2}.
\end{equation*}

\begin{itemize}

\item[(1-1)] If $S_{\min}(0)\geq0$, then $\lambda(t)$ is nondecreasing along the harmonic-Ricci flow for any $t\in[0,T]$.

\item[(1-2)] If $S_{\min}(0)>0$, then the quantity
\begin{equation*}
\left(1-S_{\min}(0)t\right)^{2\alpha}\lambda(t)
\end{equation*}
is nondecreasing along the harmonic-Ricci flow for
$T\leq\frac{1}{S_{\min}(0)}$.

\item[(1-3)] If $S_{\min}(0)<0$, then the quantity
\begin{equation*}
\left(1-S_{\min}(0)t\right)^{2\alpha}\lambda(t)
\end{equation*}
is nondecreasing along the harmonic-Ricci flow for any $t\in[0,T]$.

\end{itemize}

\item[(2)] Suppose that
\begin{equation*}
\left|\nabla_{g(t)} u(t)\right|^{2}_{g(t)}g(t)\geq2du(t)\otimes du(t).
\end{equation*}

\begin{itemize}

\item[(2-1)] If $S_{\min}(0)\geq0$, then $\lambda(t)$ is nondecreasing along the harmonic-Ricci flow for any $t\in[0,T]$.

\item[(2-2)] If $S_{\min}(0)>0$, then the quantity
\begin{equation*}
\left(1-S_{\min}(0)t\right)\lambda(t)
\end{equation*}
is nondecreasing along the harmonic-Ricci flow for
$T\leq\frac{1}{S_{\min}(0)}$.

\item[(2-3)] If $S_{\min}(0)<0$, then the quantity
\begin{equation*}
\left(1-S_{\min}(0)t\right)\lambda(t)
\end{equation*}
is nondecreasing along the harmonic-Ricci flow for any $t\in[0,T]$.

\end{itemize}
\end{itemize}
\end{corollary}

When we restrict to the Ricci flow, we obtain

\begin{corollary}\label{c1.8} Let $(g(t))_{t\in[0,T]}$ be a solution of the Ricci flow on a compact Riemannian surface $\Sigma$ and $\lambda(t)$ denote the eigenvalue of the Laplacian $\Delta_{g(t)}$.

\begin{itemize}

\item[(1)] If $R_{\min}(0)\geq0$, then $\lambda(t)$ is nondecreasing along the Ricci flow for any $t\in[0,T]$.

\item[(2)] If $R_{\min}(0)>0$, then the quantity $(1-R_{\min}(0)t)\lambda(t)$ is nondecreasing along the Ricci flow for $T\leq\frac{1}{R_{\min}(0)}$.

\item[(3)] If $R_{\min}(0)<0$, then the quantity $(1-R_{\min}(0)t)\lambda(t)$ is nondecreasing along the Ricci flow for any $t\in[0,T]$.

\end{itemize}
\end{corollary}

\begin{remark}\label{r1.9} Let $(g(t))_{t\in[0,T]}$ be a solution of the Ricci flow on a compact Riemannian surface $\Sigma$ with nonnegative scalar curvsture and $\lambda(t)$ denote the eigenvalue of the Laplaican $\Delta_{g(t)}$. Then $\lambda(t)$ is nondecreasing along the Ricci flow for any $t\in[0,T]$.
\end{remark}

Since
\begin{equation}
\mu(g,u):=\inf\left\{\mathcal{F}(g,u,f)\Big| f\in C^{\infty}(M), \ \ \int_{M}e^{-f}\!\ dV_{g}=1\right\}\label{1.11}
\end{equation}
is the smallest eigenvalue of the operator $\Delta_{g,u}:=-4\Delta_{g}+R_{g}-2\left|\nabla_{g} u\right|^{2}_{g}$, we can consider the evolution equation for this eigenvalue under the harmonic-Ricci flow. To the operator $\Delta_{g,u}$ we associate a functional
\begin{equation}
\lambda_{g,u}(f):=\int_{M}f\Delta_{g,u}f\!\ dV_{g}.\label{1.12}
\end{equation}
When $f$ is an eigenfunction of the the operator $\Delta_{g,u}$ with the eigenvalue $\lambda$ and normalized by $\int_{X}f^{2}\!\ dV_{g}=1$, we obtain $\lambda_{g,u}(f)=\lambda$. Hence we can suffice to study the evolution equation for $\frac{d}{dy}\lambda_{g,u}(f)$ under the harmonic-Ricci flow.

\begin{theorem}\label{t1.10} Suppose that $(g(t),u(t))$ is a solution of the harmonic-Ricci flow on a compact Riemannian manifold $M$ and $f(t)$ is an eigenfunction of $\Delta_{g(t),u(t)}$, i.e., $\Delta_{g(t),u(t)}f(t)=\lambda(t)f(t)$(where $\lambda(t)$
is only a function of time $t$), with the normalized condition $\int_{M}f(t)^{2}\!\ dV_{g(t)}=1$. Then we have
\begin{eqnarray}
\frac{d}{dt}\lambda(t)&=&\frac{d}{d t}\lambda_{g,u}(f(t)) \ \ = \ \ \int_{M}2\left\langle\mathcal{S}_{g(t),u(t)},df(t)\otimes
df(t)\right\rangle_{g(t)}dV_{g(t)}\nonumber\\
&&+ \ \int_{M}f(t)^{2}\left[\left|\mathcal{S}_{g(t),u(t)}\right|^{2}_{g(t)}
\label{1.13}+2\left|\Delta_{g(t)}u(t)\right|^{2}_{g(t)}
\right]dV_{g(t)}.
\end{eqnarray}
\end{theorem}

In \cite{L}, List proved the nonnegativity of the operator $\mathcal{S}_{g(t),u(t)}$ is preserved by the harmonic-Ricci flow, hence

\begin{corollary}\label{c1.11} If ${\rm Ric}_{g(0)}-2du(0)\otimes du(0)\geq0$, then the
eigenvalues of the operator $\Delta_{g(t),u(t)}$ are
nondecreasing under the harmonic-Ricci flow.
\end{corollary}

\begin{remark} \label{r1.12}If we choose $u(t)\equiv0$, then we obtain X. Cao's result \cite{C1}.
\end{remark}

There is another expression of $\frac{d}{dt}\lambda(t)$.

\begin{theorem}\label{t1.13} Suppose that $(g(t),u(t))$ is a solution of the harmonic-Ricci flow on a compact Riemannian manifold $M$ and $f(t)$ is an eigenfunction of $\Delta_{g(t),u(t)}$, i.e., $\Delta_{g(t),u(t)}f(t)=\lambda(t)f(t)$(where
$\lambda(t)$ is only a function of time $t$), with the normalized condition $\int_{M}f(t)^{2}\!\ dV_{g(t)}=1$. Then we have
\begin{eqnarray}
&&\frac{d}{dt}\lambda(t) \ \ = \ \ \frac{d}{d t}
\lambda_{g,u}(f(t)) \ \ = \ \ \frac{1}{2}\int_{M}\left|\mathcal{S}_{g(t),u(t)}+\nabla^{2}_{g(t)}\varphi(t)\right|^{2}_{g(t)}
e^{-\varphi(t)}\!\ dV_{g(t)}\nonumber\\
&&+ \ \frac{1}{4}\int_{M}\left|\mathcal{S}_{g(t),u(t)}\right|^{2}_{g(t)}
e^{-\varphi(t)}\!\ dV_{g(t)}+\int_{M}\left|\langle du(t),d\varphi(t)\rangle_{g(t)}\right|^{2}e^{-\varphi(t)}\!\ dV_{g(t)}
\nonumber\\
&&+ \ 2\int_{M}\left|\nabla^{2}_{g(t)}u(t)\right|^{2}_{g(t)}
e^{-\varphi(t)}\!\ dV_{g(t)}-\int_{M}\Delta_{g(t)}\left(\left|\nabla_{g(t)} u(t)\right|^{2}_{g(t)}\right)e^{-\varphi(t)}\!\ dV_{g(t)}\label{1.14}\\
&&+\frac{1}{4}\int_{M}\left|\mathcal{S}_{g(t),u(t)}+4du(t)\otimes du(t)\right|^{2}_{g(t)}e^{-\varphi(t)}\!\ dV_{g(t)},\nonumber
\end{eqnarray}
where $f(t)^{2}=e^{-\varphi(t)}$.
\end{theorem}

\begin{remark}\label{r1.14} When $u\equiv0$, (\ref{1.14}) reduces to
J. Li's formula \cite{Li}.
\end{remark}

Suppose that $M$ is a compact manifold of dimension $n$. For any Riemannian metric $g$, any smooth functions $u,f$, and any positive number $\tau$, we define
\begin{equation}
\mathcal{W}_{\pm}(g,u,f,\tau):=\int_{M}\left[\tau\left(S_{g}
+\left|\nabla_{g} f\right|^{2}_{g}\right)\mp f\pm n\right]\frac{e^{-f}}{(4\pi\tau)^{n/2}}\!\ dV_{g}.\label{1.15}
\end{equation}
Set
\begin{eqnarray*}
\mu_{\pm}(g,u,\tau)&:=&\inf\left\{\mathcal{W}_{\pm}(g,u,f,\tau)\Big| f\in C^{\infty}(M), \ \ \ \int_{M}\frac{e^{-f}}{(4\pi\tau)^{n/2}}\!\ dV_{g}=1\right\},\\
\nu_{-}(g,u)&:=&\inf\{\mu_{-}(g,u,\tau)|\tau>0\}, \ \ \
\nu_{+}(g,u) \ \ := \ \ \sup\{\mu_{+}(g,u,\tau)|\tau>0\}.
\end{eqnarray*}

The first variation of $\nu_{\pm}(g(s),u(s))$ is

\begin{theorem}\label{t1.15} Suppose that $(M,g)$ is a compact Riemannian manifold and $u$ a smooth function on $M$. Let $h$ be any
symmetric covariant $2$-tensor on $M$ and set $g(s):=g+sh$. Let $v$ be any smooth function on $M$ and $u(s):=u+sv$. If $\nu_{\pm}(g(s),u(s))=\mathcal{W}_{\pm}(g(s),u(s),f_{\pm}(s),\tau_{\pm}(s))$ for some smooth functions $f_{\pm}(s)$ with $\int_{M}e^{-f_{\pm}(s)}dV_{g}/(4\pi\tau_{\pm}(s))^{n/2}=1$ and constants $\tau_{\pm}(s)>0$, then
\begin{eqnarray}
\frac{d}{ds}\Big|_{s=0}\nu_{\pm}(g(s),u(s))
&=&4\tau_{\pm}\int_{M}v\left(\Delta_{g} u-\langle du,df_{\pm}\rangle_{g}\right)\frac{e^{-f_{\pm}}}{
(4\pi\tau_{\pm})^{n/2}}dV_{g}\nonumber\\
&&- \ \tau_{\pm}\int_{M}
\left(\left\langle h,\mathcal{S}_{g,u}\right\rangle_{g}+\left\langle h,\nabla^{2}_{g}f\right\rangle_{g}\pm\frac{{\rm tr}_{g}h}{2\tau_{\pm}}\right)\frac{e^{-f_{\pm}}dV_{g}}{
(4\pi\tau_{\pm})^{n/2}},\label{1.16}
\end{eqnarray}
where $f_{\pm}:=f_{\pm}(0)$ and $\tau_{\pm}:=\tau_{\pm}(0)$. In particular, the critical points of $\nu_{\pm}(\cdot,\cdot)$ satisfy
\begin{equation*}
\mathcal{S}_{g,u}+\nabla^{2}_{g}f\pm\frac{1}{2\tau_{\pm}}g=0, \ \ \ \Delta_{g}u=\left\langle du,df_{\pm}\right\rangle_{g}.
\end{equation*}
Consequently, if $\mathcal{W}_{\pm}(g,u,f,\tau)$ and
$\nu_{\pm}(g,u)$ achieve their extremum, then $(M,g)$ is a gradient expanding and shrinker harmonic-Ricci soliton according to the sign.
\end{theorem}

\begin{corollary}\label{c1.16} Suppose that $(M,g)$ is a compact Riemannian
manifold and $u$ a smooth function on $M$. Let $h$ be any symmetric covariant $2$-tensor on $M$ and set $g(s):=g+sh$. Let $v$ be any smooth function on $M$ and $u(s):=u+sv$. If $\nu_{\pm}(g(s),u(s))=\mathcal{W}_{\pm}(g(s),u(s),f_{\pm}(s),\tau_{\pm}(s))$ for some smooth function $f_{\pm}(s)$ with $\int_{M}e^{-f_{\pm}(s)}dV/(4\pi\tau_{\pm}(s))^{n/2}=1$ and a constant
$\tau_{\pm}(s)>0$, and $(g,u)$ is a critical point of
$\nu_{\pm}(\cdot,\cdot)$, then
\begin{equation*}
{\rm Ric}_{g}=\mp\frac{1}{2\tau_{\pm}}g, \ \ \ f_{\pm}\equiv{\rm constant}, \ \ \ u\equiv{\rm constant}.
\end{equation*}
Thus, if $\mathcal{W}_{\pm}(g,u,\cdot,\cdot)$ achieve their minimum and $(g,u)$ is
a critical point of $\nu_{\pm}(\cdot,\cdot)$, then $(M,g)$ is an Einstein manifold and $u$ is a constant function.
\end{corollary}

\begin{remark}\label{r1.17} In the situation of Corollary \ref{c1.16}, by normalization, we my choose $f_{\pm}=\frac{n}{2}$ and $u=0$.
\end{remark}

{\bf Acknowledgements.} Part of work was done when the author visited Center of Mathematical Science at Zhejiang University in 2010. The author would like to thank Professor Kefeng Liu, who teaches the author mathematics. Furthermore, I also thank Professor Hongwei Xu and other staffs in Center of Mathematical Science.

\section{Notation and commuting identities}\label{section2}

Let $M$ be a compact Riemannian manifold of dimension $n$. For any vector bundle $E$ over $M$, we denote by $\Gamma(M,E)$ the space of smooth sections of $E$. Set
\begin{eqnarray*}
\odot^{2}(M)&:=&\{v=(v_{ij})\in\Gamma(M,T^{\ast}M\otimes T^{\ast}M)|v_{ij}=v_{ji}\},\\
\odot^{2}_{+}(M)&:=&\{g=(g_{ij})\in\odot^{2}(M)|g_{ij}>0\}.
\end{eqnarray*}
Thus, $\odot^{2}(M)$ is the space of all symmetric covariant $2$-tensors on $M$ while $\odot^{2}_{+}(M)$ the space of all Riemannian metrics on $M$. The space of all smooth functions on $M$ is denoted by $C^{\infty}(M)$.

For a given Riemannian metric $g\in\odot^{2}_{+}(M)$, the corresponding Levi-civita connection $\Gamma_{g}=(\Gamma^{k}_{ij})$ is given by
\begin{equation}
\Gamma^{k}_{ij}=\frac{1}{2}g^{k\ell}\left(\partial_{i}g_{j\ell}
+\partial_{j}g_{il}-\partial_{\ell}g_{ij}\right)\label{2.1}
\end{equation}
where $\partial_{i}:=\frac{\partial}{\partial x^{i}}$ for a local coordinate system $\{x^{1},\cdots, x^{n}\}$. The Riemann tensor ${\rm Rm}_{g}=(R^{k}_{ijl})$ is determined by
\begin{equation}
R^{k}_{ij\ell}=\partial_{i}\Gamma^{k}_{j\ell}
-\partial_{j}\Gamma^{k}_{i\ell}
+\Gamma^{k}_{ip}\Gamma^{p}_{j\ell}-\Gamma^{k}_{jp}\Gamma^{p}_{i\ell}.\label{2.2}
\end{equation}
The Ricci curvature ${\rm Ric}_{g}=(R_{ij})$ is
\begin{equation}
R_{ij}=g^{k\ell}R^{\ell}_{kij}.\label{2.3}
\end{equation}
The scalar curvature $R_{g}$ of the metric $g$ now is given by
\begin{equation}
R_{g}=g^{ij}R_{ij}.\label{2.4}
\end{equation}

For any tensor $A=(A^{k_{1}\cdots k_{q}}_{j_{1}\cdots j_{p}})$ the covariant
derivative of $A$ is
\begin{equation*}
\nabla_{i}A^{k_{1}\cdots k_{q}}_{j_{1}\cdots j_{p}}=\partial_{i}A^{k_{1}\cdots k_{q}}_{j_{1}\cdots j_{p}}-\sum^{p}_{r=1}\Gamma^{m}_{ij_{r}}A^{k_{1}\cdots k_{q}}_{j_{1}\cdots m\cdots j_{p}}+\sum^{q}_{s=1}\Gamma^{k_{s}}_{im}A^{k_{1}\cdots m\cdots k_{q}}_{j_{1}\cdots j_{p}}.
\end{equation*}
Next we recall the Ricci identity:
\begin{equation*}
\nabla_{i}\nabla_{j}A^{\ell_{1}\cdots \ell_{q}}_{k_{1}\cdots k_{p}}-\nabla_{j}\nabla_{i}A^{\ell_{1}\cdots \ell_{q}}_{k_{1}\cdots k_{p}}
=\sum^{q}_{r=1}R^{l_{r}}_{ijm}A^{\ell_{1}\cdots m\cdots \ell_{q}}_{k_{1}\cdots k_{p}}
-\sum^{p}_{s=1}R^{m}_{ijk_{s}}A^{\ell_{1}\cdots \ell_{q}}_{k_{1}\cdots m\cdots k_{p}}.
\end{equation*}
In particular, for any smooth function $f\in C^{\infty}(M)$ we have
\begin{equation*}
\nabla_{i}\nabla_{j}f=\nabla_{j}\nabla_{i}f.
\end{equation*}
The Bianchi identities are
\begin{eqnarray}
0&=&R_{ijk\ell}+R_{iklj}+R_{i\ell jk},\label{2.5}\\
0&=&\nabla_{q}R_{ijk\ell}+\nabla_{i}R_{jqk\ell}
+\nabla_{j}R_{qik\ell}\label{2.6}
\end{eqnarray}
and the contracted Bianchi identities are
\begin{eqnarray}
0&=&2\nabla^{j}R_{ij}-\nabla_{i}R_{g},\label{2.7}\\
0&=&\nabla_{i}R_{jk}-\nabla_{j}R_{ik}
+\nabla^{\ell}R_{\ell kij}.\label{2.8}
\end{eqnarray}

\section{Harmonic-Ricci flow and the evolution equations}
\label{section3}

Motivated by static Einstein vacuum equation, List \cite{L} introduced the harmonic-Ricci flow(Originally, it is called the Ricci flow coupled with the harmonic map flow.). Such a flow is similar to the Ricci flow and is the following coupled system
\begin{eqnarray}
\frac{\partial}{\partial t}g(t)&=&-2{\rm Ric}_{g(t)}
+4du(t)\otimes du(t), \label{3.1}\\
\frac{\partial}{\partial t}u(t)&=&\Delta_{g(t)}u(t)\label{3.2}
\end{eqnarray}
for a family of Riemannian metrics $g(t)$ and a family of smooth functions $u(t)$. Locally, we have
\begin{equation}
\frac{\partial}{\partial t}g_{ij}=-2R_{ij}+4\partial_{i}u\cdot\partial_{j}u, \ \ \ \frac{\partial}{\partial t}u=\Delta_{g(t)}u(t).\label{3.3}
\end{equation}
Introduce a new symmetric tensor field $\mathcal{S}_{g(t),u(t)}=(S_{ij})\in\odot^{2}(M)$ by
\begin{equation}
S_{ij}:=R_{ij}-2\partial_{i}u\cdot\partial_{j}u.\label{3.4}
\end{equation}
Then its trace $S_{g(t),u(t)}$ is equal to
\begin{equation}
S_{g(t),u(t)}=g^{ij}S_{ij}=R_{g(t)}-2\left|\nabla_{g(t)} u(t)\right|^{2}_{g(t)}.\label{3.5}
\end{equation}

The evolution equation for $R_{g(t)}$ is
\begin{eqnarray}
\frac{\partial}{\partial t}R_{g(t)}&=&\Delta_{g(t)}R_{g(t)}+2|{\rm Ric}_{g(t)}|^{2}_{g(t)}\nonumber\\
&&+4\left|\Delta_{g(t)}u(t)\right|^{2}_{g(t)}-4\left|\nabla^{2}_{g(t)}u(t)\right|^{2}_{g(t)}
-8\left\langle{\rm Ric}_{g(t)},du(t)\otimes du(t)\right\rangle_{g(t)}.\label{3.6}
\end{eqnarray}
Also, we have the evolution equation for $\left|\nabla_{g(t)} u\right|^{2}_{g(t)}$:
\begin{equation}
\frac{\partial}{\partial t}\left|\nabla_{g(t)} u(t)\right|^{2}_{g(t)}=\Delta_{g(t)}\left|\nabla_{g(t)} u(t)\right|^{2}_{g(t)}-2\left|\nabla^{2}_{g(t)}u(t)
\right|^{2}_{g(t)}-4\left|\nabla_{g(t)} u(t)\right|^{4}_{g(t)},\label{3.7}
\end{equation}
and the evolution equation for $S_{g(t),u(t)}$:
\begin{equation}
\frac{\partial}{\partial t}S_{g(t),u(t)}=\Delta_{g(t)}S_{g(t),u(t)}+2\left|\mathcal{S}_{g(t),u(t)}\right|^{2}_{g(t)}+
4\left|\Delta_{g(t)}u(t)\right|^{2}_{g(t)}.\label{3.8}
\end{equation}

\section{Entropies for harmonic-Ricci flow}\label{section4}

Motivated by Perelman's entropy, List \cite{L} introduced the similar functional for the harmonic-Ricci flow:
\begin{equation*}
\odot^{2}_{+}(M)\times C^{\infty}(M)\times C^{\infty}(M)\longrightarrow
\mathbb{R},\ \ \ (g,u,f)\longmapsto\mathcal{F}(g,u,f)
\end{equation*}
where
\begin{equation}
\mathcal{F}(g,u,f):=\int_{M}\left(R_{g}+\left|\nabla_{g} f\right|^{2}_{g}-2\left|\nabla_{g} u\right|^{2}_{g}\right)e^{-f}\!\ dV_{g}.\label{4.1}
\end{equation}
He also showed that if $(g(t),u(t),f(t))$ satisfies the following system
\begin{eqnarray}
\frac{\partial}{\partial t}g(t)&=&-2{\rm Ric}_{g(t)}+4du(t)\otimes du(t)-2\nabla^{2}_{g(t)}f(t),\nonumber\\
\frac{\partial}{\partial t}u(t)&=&\Delta_{g(t)}u(t)-\langle du(t),df(t)\rangle_{g(t)},\label{4.2}\\
\frac{\partial}{\partial t}f(t)&=&-\Delta_{g(t)}f(t)-R_{g(t)}+2\left|\nabla_{g(t)} u(t)\right|^{2}_{g(t)},\nonumber
\end{eqnarray}
then the evolution of the entropy is given by
\begin{eqnarray}
\frac{d}{dt}\mathcal{F}(g(t),u(t),f(t))&=&2\int_{M}
\bigg(\left|\mathcal{S}_{g(t),u(t)}+\nabla^{2}_{g(t)}
f(t)\right|^{2}_{g(t)}\nonumber\\
&&+ \ 2\left|\Delta_{g(t)}u(t)-\langle du(t),df(t)\rangle_{g(t)}\right|^{2}_{g(t)}\bigg)e^{-f(t)}dV_{g(t)}
\geq0.\label{4.3}
\end{eqnarray}

\begin{remark}\label{r4.1} The above system (\ref{4.2}) is equivalent to the following
\begin{eqnarray}
\frac{\partial}{\partial t}g(t)&=&-2{\rm Ric}_{g(t)}+4du(t)\otimes du(t),\nonumber\\
\frac{\partial}{\partial t}u(t)&=&\Delta_{g(t)}u(t),\label{4.4}\\
\frac{\partial}{\partial t}f(t)&=&-\Delta_{g(t)}f(t)-R_{g(t)}
+\left|\nabla_{g(t)} f(t)\right|^{2}_{g(t)}+
2\left|\nabla_{g(t)} u(t)\right|^{2}_{g(t)}.\nonumber
\end{eqnarray}
The same evolution of the entropy holds for this system (\ref{4.4}).
\end{remark}

In particular, the entropy is nondecreasing and the equality holds if and only if $(g(t),u(t),f(t))$ satisfies
\begin{equation}
\mathcal{S}_{g(t),u(t)}+\nabla^{2}_{g(t)}f(t)=0, \ \ \ \Delta_{g(t)}u(t)-\langle du(t),df(t)\rangle_{g(t)}=0.\label{4.5}
\end{equation}

\begin{definition}\label{d4.2} The $\mathcal{E}$-functional is defined as
\begin{equation*}
\odot^{2}_{+}(M)\times C^{\infty}(M)\times C^{\infty}(M)\longrightarrow\mathbb{R}, \ \ \ (g,u,f)\longmapsto\mathcal{E}(g,u,f),
\end{equation*}
where
\begin{equation}
\mathcal{E}(g,u,f):=\int_{M}\left(R_{g}-2\left|\nabla_{g} u\right|^{2}_{g}\right)e^{-f}\!\ dV_{g}.\label{4.6}
\end{equation}
\end{definition}

\begin{proposition}\label{p4.3} Under the evolution equation (\ref{4.4}), one has
\begin{eqnarray}
\frac{d}{dt}\mathcal{E}(g(t),u(t),f(t))&=&
2\int_{M}\left|\mathcal{S}_{g(t),u(t)}\right|^{2}_{g(t)}e^{-f(t)}\!\ dV_{g(t)}\label{4.7}\\
&&+ \ 4\int_{M}\left|\Delta_{g(t)}u(t)\right|^{2}_{g(t)}e^{-f(t)}\!\ dV_{g(t)}.\nonumber
\end{eqnarray}
\end{proposition}

\begin{proof} Since $S_{g(t),u(t)}=R_{g(t)}
-2\left|\nabla_{g(t)} u(t)\right|^{2}_{g(t)}$ and
\begin{eqnarray*}
\frac{\partial}{\partial t}S_{g(t),u(t)}&=&\Delta_{g(t)}S_{g(t),u(t)}+2\left|\mathcal{S}_{g(t),u(t)}\right|^{2}_{g(t)}+4\left|\Delta_{g(t)}u(t)\right|^{2}_{g(t)},\\
\frac{\partial}{\partial t}dV_{g(t)}&=&-S_{g(t),u(t)}\!\ dV_{g(t)},
\end{eqnarray*}
we have
\begin{eqnarray*}
&&\frac{d}{dt}\mathcal{E}(g(t),u(t),f(t))\\
&=&\int_{M}\left(\frac{\partial}{\partial t}S_{g(t),u(t)}\right)e^{-f(t)}
\!\ dV_{g(t)}+\int_{M}S_{g(t),u(t)}\frac{\partial}{\partial t}\left(e^{-f(t)}\!\ dV_{g(t)}\right)\\
&=&\int_{M}\left(\Delta_{g(t)}S_{g(t),u(t)}+2\left|\mathcal{S}_{g(t),u(t)}\right|^{2}_{g(t)}
+4\left|\Delta_{g(t)}u(t)\right|^{2}_{g(t)}\right)e^{-f(t)}\!\ dV_{g(t)}\\
&&+ \ \int_{M}S_{g(t),u(t)}\left(-\frac{\partial}{\partial t}f(t)-S_{g(t),u(t)}
\right)e^{-f(t)}\!\ dV_{g(t)}\\
&=&2\int_{M}\left|\mathcal{S}_{g(t),u(t)}\right|^{2}_{g(t)}
e^{-f(t)}\!\ dV_{g(t)}
+4\int_{M}\left|\Delta_{g(t)}u(t)\right|^{2}_{g(t)}e^{-f(t)}\!\ dV_{g(t)}\\
&&- \ \int_{M}S_{g(t),u(t)}\bigg(\Delta_{g(t)}f(t)-\left|
\nabla_{g(t)} f(t)\right|^{2}_{g(t)}+\frac{\partial}{\partial t}f(t)+S_{g(t),u(t)}\bigg)
e^{-f(t)}dV_{g(t)}
\end{eqnarray*}
which implies (\ref{4.7}).
\end{proof}

\begin{definition} \label{d4.4}For any $k\geq1$ we define
\begin{equation}
\mathcal{F}_{k}(g,u,f):=\int_{M}\left(kR_{g}+\left|\nabla_{g} f\right|^{2}_{g}-2k\left|\nabla_{g} u\right|^{2}_{g}\right)
e^{-f}\!\ dV_{g}.\label{4.8}
\end{equation}
By definition, it is easy to show that
\begin{equation}
\mathcal{F}_{k}(g,u,f)=(k-1)\mathcal{E}(g,u,f)+\mathcal{F}(g,u,f).\label{4.9}
\end{equation}
When $k=1$, this is the $\mathcal{F}$-functional.
\end{definition}

\begin{theorem}\label{t4.5} Under the evolution equation (\ref{4.4}), one has
\begin{eqnarray}
\frac{d}{dt}\mathcal{F}_{k}(g(t),u(t),f(t))&=&2(k-1)\int_{M}
\left|\mathcal{S}_{g(t),u(t)}\right|^{2}_{g(t)}
e^{-f(t)}dV_{g(t)}\nonumber\\
&&+ \ 2\int_{M}\left|\mathcal{S}_{g(t),u(t)}+\nabla^{2}_{g(t)}f(t)\right|^{2}_{g(t)}e^{-f(t)}dV_{g(t)}\label{4.10}\\
&&+ \ 4(k-1)\int_{M}\left|\Delta_{g(t)}u(t)\right|^{2}_{g(t)}e^{-f(t)}dV_{g(t)}\nonumber\\
&&+ \ 4\int_{M}\left|\Delta_{g(t)}u(t)-\langle du(t),df(t)\rangle_{g(t)}
\right|^{2}_{g(t)}e^{-f(t)}dV_{g(t)}.\nonumber
\end{eqnarray}
Furthermore, the monotonicity is strict unless $g(t)$ is Ricci-flat, $u(t)$ is constant and $f(t)$ is constant.
\end{theorem}

\begin{proof} It immediately follows from (\ref{4.3}) and (\ref{4.7}).
\end{proof}

Set
\begin{equation}
\mu_{k}(g,u):=\inf\left\{\mathcal{F}_{k}(g,u,f)\Big|f\in C^{\infty}(M), \ \ \int_{M}e^{-f}dV_{g}=1\right\}.\label{4.11}
\end{equation}
Then $\mu_{k}(g,u)$ is the lowest eigenvalue of $-4\Delta_{g}+k\left(R_{g}
-2\left|\nabla_{g} u\right|^{2}_{g}\right)$.

\section{Compact steady harmonic-Ricci breathers}
\label{section5}

In this section we give an alternative proof on some results on compact steady harmonic-Ricci breathers that were proved in \cite{L,M2}.

\begin{definition} \label{d5.1}A solution $(g(t),u(t))$ of the harmonic-Ricci flow (\ref{1.1})--(\ref{1.2}) is called a
{\it harmonic-Ricci breather} if there exist $t_{1}<t_{2}$, a diffeomorphism $\psi: M\to M$ and a constant $\alpha>0$ such that
\begin{equation*}
g(t_{2})=\alpha\psi^{\ast}g(t_{1}), \ \ \ u(t_{2})=\psi^{\ast}u(t_{1}).
\end{equation*}
The case $\alpha<1,\alpha=1$, and $\alpha>1$, correspond to {\it
shrinking}, {\it steady} and {\it expanding harmonic-Ricci breathers}.
\end{definition}

\begin{theorem}\label{t5.2} If $(g(t),u(t))$ is a solution of the harmonic-Ricci flow on a compact Riemannian manifold $M$, then the lowest eigenvalue
$\mu_{k}(g(t),u(t))$ of the operator $-4\Delta_{g(t)}+k(R_{g(t)}-2|\nabla_{g(t)} u(t)|^{2}_{g(t)})$ is nondecreasing under the harmonic-Ricci flow. The monotonicity is streat unless $g(t)$ is Ricci-flat and $u(t)$ is constant.
\end{theorem}

\begin{proof} The proof is similar to that given in \cite{Li}. For any $t_{1}<t_{2}$, suppose that
\begin{equation*}
\mu_{k}(g(t_{2}),u(t_{2}))=\mathcal{F}_{k}(g(t_{2}),u(t_{2}),f_{k}(t_{2}))
\end{equation*}
for some smooth function $f_{k}(x)$. Being an initial value, $f_{k}(x)
=f_{k}(x,t_{2})$ for some smooth function $f_{k}(x,t)$ satisfying the evolution equation (\ref{4.4}). The monotonicity formula (\ref{4.10}) implies
\begin{equation*}
\mu_{k}(g(t_{2}),u(t_{2}))\geq
\mathcal{F}_{k}(g(t_{1}),u(t_{1}),f_{k}(t_{1}))\geq\mu_{k}(g(t_{1}),u(t_{1})).
\end{equation*}
This completes the proof.
\end{proof}

\begin{corollary} \label{c5.3}On a compact Riemannian manifold, the
lowest eigenvalues of $-\Delta_{g(t)}+\frac{1}{2}(R_{g(t)}-2|\nabla_{g(t)} u(t)|^{2}_{g(t)})$ are nondecreasing under the harmonic-Ricci flow.
\end{corollary}

\begin{proof} Since $\mu_{2}(g(t),u(t))/4$ is the lowest eigenvalue of the above operator, the result immediately follows from Theorem \ref{t5.2}.
\end{proof}

\begin{corollary} \label{c5.4}There is no compact steady harmonic-Ricci breather other
than $(M,g(t))$ is Ricci-flat and $u$ is constant.
\end{corollary}

\begin{proof} If $(g(t),u(t))$ is a steady harmonic-Ricci breather, then for $t_{1}<t_{2}$ given in the definition, we have
\begin{equation*}
\mu_{k}(g(t_{1}),u(t_{1}))=\mu_{k}(g(t_{2}),u(t_{2}))
\end{equation*}
hence, using Theorem \ref{t5.2}, for any $t\in[t_{1},t_{2}]$ we must have
\begin{equation*}
\frac{d}{dt}\mu_{k}(g(t),u(t))\equiv0.
\end{equation*}
Thus $(M,g(t))$ is Ricci-flat and $u(t)$ is constant.
\end{proof}

\section{Compact expanding harmonic-Ricci breathers}
\label{section6}

Inspired by \cite{Li}, we define a new functional
\begin{equation*}
\odot^{2}_{+}(M)\times C^{\infty}(M)\times C^{\infty}(\mathbb{R})\times
C^{\infty}(M)
\longrightarrow\mathbb{R}, \ \ \ (g,u,\tau,f)\longmapsto\mathcal{W}
_{+}
(g,u,\tau,f),
\end{equation*}
where ($\tau=\tau(t), t\in\mathbb{R}$)
\begin{equation}
\mathcal{W}_{+}(g,u,\tau,f):=\tau^{2}\int_{M}\left(R_{g}
+\frac{n}{2\tau}+\Delta_{g}f-2\left|\nabla_{g} u\right|^{2}_{g}\right)e^{-f}dV_{g}.\label{6.1}
\end{equation}
Similarly, we define a family of functionals
\begin{equation}
\mathcal{W}_{+,k}(g,u,\tau,f):=\tau^{2}\int_{M}\left[k
\left(R_{g}+\frac{n}{2\tau}\right)+\Delta_{g}f-2k\left|\nabla_{g}
u\right|^{2}_{g}\right]e^{-f}dV_{g}.\label{6.2}
\end{equation}
It's clear that $\mathcal{W}_{+,1}(g,u,\tau,f)=
\mathcal{W}_{+}(g,u,\tau,f)$.

\begin{lemma} \label{l6.1}One has
\begin{eqnarray*}
\mathcal{W}_{+}(g,u,\tau,f)&=&\tau^{2}\mathcal{F}(g,u,f)
+\frac{n}{2}\tau\int_{M}e^{-f}dV_{g},\\
\mathcal{W}_{+,k}(g,u,\tau,f)&=&\tau^{2}\mathcal{F}_{k}(g,u,f)
+\frac{kn}{2}\tau\int_{M}e^{-f}dV_{g},\\
\mathcal{W}_{+,k}(g,u,\tau,f)&=&
\mathcal{W}_{+}(g,u,\tau,f)\\
&&+(k-1)\left(\tau^{2}\mathcal{E}(g,u,f)
+\frac{n}{2}\tau\int_{M}e^{-f}dV_{g}\right).
\end{eqnarray*}
\end{lemma}

\begin{proof} Since $\Delta(e^{-f})=(-\Delta f+|\nabla f|^{2})e^{-f}$, it follows that
\begin{eqnarray*}
&&\mathcal{W}_{+}(g,u,\tau,f)-\tau^{2}\mathcal{F}(g,u,f)\\
&=&\frac{n}{2}\tau\int_{M}e^{-f}dV_{g}+\tau^{2}\int_{M}\left(\Delta_{g}f
-\left|\nabla_{g} f\right|^{2}_{g}\right)e^{-f}dV_{g}\\
&=&\frac{n}{2}\tau\int_{M}e^{-f}dV_{g}+\tau^{2}\int_{M}\Delta_{g}
\left(e^{-f}\right)dV_{g} \ \ = \ \ \frac{n}{2}\tau\int_{M}e^{-f}dV_{g}.
\end{eqnarray*}
Similarly, we can prove the rest two relations.
\end{proof}

\begin{theorem} \label{t6.2}Under the following coupled system
\begin{eqnarray*}
\frac{\partial }{\partial t}g(t)&=&-2{\rm Ric}_{g(t)}+4du(t)\otimes
du(t)-2\nabla^{2}_{g(t)}f(t),\\
\frac{\partial}{\partial t}u(t)&=&\Delta_{g(t)}u(t)-\langle du(t),df(t)\rangle_{g(t)},\\
\frac{\partial}{\partial t}f(t)&=&-\Delta_{g(t)}f(t)-R_{g(t)}+2\left|\nabla_{g(t)} u(t)\right|^{2}_{g(t)},\\
\frac{d}{dt}\tau(t)&=&1,
\end{eqnarray*}
the first variation formula for $\mathcal{W}_{+}(g(t),u(t),\tau(t),f(t))$ is
\begin{eqnarray}
&&\frac{d}{dt}\mathcal{W}_{+}(g(t),u(t),\tau(t),f(t))\label{6.3}\\
&=&2\tau(t)^{2}\int_{M}\left|\mathcal{S}_{g(t),u(t)}
+\nabla^{2}_{g(t)}f(t)+\frac{1}{2\tau(t)}
g(t)\right|^{2}_{g(t)}e^{-f(t)}dV_{g(t)}\nonumber\\
&&+ \ 4\tau(t)^{2}\int_{M}\left|\Delta_{g(t)}u(t)-\langle du(t),df(t)\rangle_{g(t)}\right|^{2}_{g(t)}e^{-f(t)}dV_{g(t)},\nonumber
\end{eqnarray}
and the first variation formula for $\mathcal{W}_{+,k}(g(t),u(t),\tau(t),f(t))$ is
\begin{eqnarray}
&&\frac{d}{dt}\mathcal{W}_{+,k}(g(t),u(t),\tau(t),f(t))\label{6.4}\\
&=&2\tau(t)^{2}\int_{M}\left|\mathcal{S}_{g(t),u(t)}
+\nabla^{2}_{g(t)}f(t)+\frac{1}{2\tau(t)}
g(t)\right|^{2}_{g(t)}e^{-f(t)}dV_{g(t)}\nonumber\\
&&+2(k-1)\tau(t)^{2}\int_{M}\left|\mathcal{S}_{g(t),u(t)}+\frac{1}{2\tau(t)}g(t)\right|^{2}_{g(t)}e^{-f(t)}dV_{g(t)}\nonumber\\
&&+4\tau(t)^{2}\int_{M}\left|\Delta_{g(t)}u(t)-\langle du(t),df(t)\rangle_{g(t)}\right|^{2}_{g(t)}e^{-f(t)}dV_{g(t)}\nonumber\\
&&+4(k-1)\tau(t)^{2}\int_{M}\left|\Delta_{g(t)}u(t)\right|^{2}_{g(t)}e^{-f(t)}dV_{g(t)}.\nonumber
\end{eqnarray}

\end{theorem}

\begin{proof} Under the above coupled system, we first observe that
\begin{equation}
\frac{d}{dt}\left(\int_{M}e^{-f(t)}dV_{g(t)}\right)=0.\label{6.5}
\end{equation}
In fact, from $\frac{\partial}{\partial t}dV_{g(t)}=[-S_{g(t),u(t)}
-\Delta_{g(t)}f(t)dV_{g(t)}$ we obtain
\begin{eqnarray*}
\frac{d}{dt}\left(\int_{M}e^{-f(t)}dV_{g(t)}\right)&=&\int_{M}\left(-\frac{\partial}{\partial t}f(t)\cdot dV_{g(t)}
+\frac{\partial}{\partial t}dV_{g(t)}\right)e^{-f(t)}\\
&=&\int_{M}\left[\Delta_{g(t)}f(t)+S_{g(t),u(t)}\right.\\
&& \ -\left.S_{g(t),u(t)}-
\Delta_{g(t)}f(t)\right]e^{-f(t)}dV_{g(t)} \\
&=&0.
\end{eqnarray*}
Lemma \ref{l6.1} and the identity (\ref{6.5}) implies
\begin{eqnarray*}
&&\frac{d}{dt}\mathcal{W}_{+}(g(t),u(t),\tau(t),f(t))\\
&=&\tau(t)^{2}\frac{d}{dt}\mathcal{F}(g(t),u(t),f(t))+2\tau(t)\mathcal{F}
(g(t),u(t),f(t))+\frac{n}{2}\int_{M}e^{-f(t)}dV_{g(t)}\\
&=&2\tau(t)^{2}\int_{M}\left|\mathcal{S}_{g(t),u(t)}
+\nabla^{2}_{g(t)}f(t)\right|^{2}_{g(t)}e^{-f(t)}
dV_{g(t)}\\
&&+ \ 4\tau(t)^{2}\int_{M}\left|\Delta_{g(t)}u(t)-\langle du(t),df(t)\rangle_{g(t)}\right|^{2}e^{-f(t)}dV_{g(t)}\\
&&+ \ 2\tau(t)\int_{M}\left(S_{g(t),u(t)}
+\left|\nabla_{g(t)} f(t)\right|^{2}_{g(t)}\right)e^{-f(t)}dV_{g(t)}+\frac{n}{2}\int_{M}e^{-f(t)}dV_{g(t)}
\end{eqnarray*}
which is (\ref{6.3}). Using Lemma \ref{l6.1} and the same method we
can prove (\ref{6.4}).
\end{proof}

\begin{remark} \label{r6.3}Under the following coupled system
\begin{eqnarray*}
\frac{\partial}{\partial t}g(t)&=&-2{\rm Ric}_{g(t)}+4du(t)\otimes du(t),\\
\frac{\partial}{\partial t}u(t)&=&\Delta_{g(t)}u(t),\\
\frac{\partial}{\partial t}f(t)&=&-\Delta_{g(t)}f(t)+\left|\nabla_{g(t)} f(t)\right|^{2}_{g(t)}-R_{g(t)}+2\left|\nabla_{g(t)} u(t)\right|^{2}_{g(t)},\\
\frac{d}{dt}\tau(t)&=&1,
\end{eqnarray*}
the same formulas (\ref{6.3}) and (\ref{6.4}) hold for
$\mathcal{W}_{+}$ and $\mathcal{W}_{+,k}$.
\end{remark}

Define
\begin{equation}
\mu_{+}(g,u,\tau):=\inf\left\{\mathcal{W}_{+}(g,u,\tau,f)\Big| f\in C^{\infty}(M), \ \ \int_{M}e^{-f}dV_{g}=1\right\}.\label{6.6}
\end{equation}

\begin{lemma} \label{l6.4}For any $\alpha>0$, one has
\begin{equation}
\mu_{+}(\alpha g,u,\alpha\tau)=\alpha\mu_{+}(g,u,\tau).\label{6.7}
\end{equation}
\end{lemma}

\begin{proof} If we set $\bar{g}:=\alpha g$, then $R_{\bar{g}}=\alpha^{-1}R_{g}$, $\Delta_{\bar{g}}f=\alpha^{-1}\Delta_{g}f$, and $|\nabla_{\bar{g}} u|^{2}_{\bar{g}}=\alpha^{-1}
|\nabla_{g(t)} u|^{2}_{g}$. Hence
\begin{eqnarray*}
\mathcal{W}_{+}(\bar{g},u,\alpha\tau,f)&=&\alpha^{2}\tau^{2}\int_{M}\left(R_{\overline{g}}+\frac{n}{2\alpha\tau}+
\Delta_{\bar{g}}f-2\left|
\nabla_{\bar{g}} u\right|^{2}_{\bar{g}}\right)e^{-f}dV_{\bar{g}}\\
&=&\alpha\tau^{2}\int_{M}\left(R_{g}+\frac{n}{2\tau}
+\Delta_{g}f-2\left|\nabla_{g(t)} u\right|^{2}_{g}\right)\alpha^{n/2}e^{-f}dV_{g}.
\end{eqnarray*}
Since $f\mapsto f-\frac{n}{2}\ln\alpha$ is one-to-one and onto, by taking the infimum we derive $\mu_{+}(\alpha g,u,\alpha\tau)=\alpha\mu_{+}(g,u,\tau)$.
\end{proof}

\begin{definition} \label{d6.5}A solution $(g(t),u(t))$ of the harmonic-Ricci flow is called a {\it harmonic-Ricci soliton} if there exists an one-parameter family of diffeomorphisms $\psi_{t}: M\to M$, satisfying $\psi_{0}={\rm id}_{M}$, and a positive scaling function $\alpha(t)$ such that
\begin{equation*}
g(t)=\alpha(t)\psi^{\ast}_{t}g(0), \ \ \ u(t)=\psi^{\ast}_{t}u(0).
\end{equation*}
The case $\frac{\partial}{\partial t}\alpha(t)=\dot{\alpha}<0$,
$\dot{\alpha}=0$, and $\dot{\alpha}>0$ correspond to
{\it shrinking}, {\it steady}, and {\it expanding harmonic-Ricci solitons},
respectively. If the diffeomorphisms $\psi_{t}$ are generated by a
(possibly time-dependent) vector field $X(t)$ that is the gradient of some function $f(t)$ on $M$, then the soliton is called {\it gradient harmonic-Ricci soliton} and $f$ is called the {\it potential of the harmonic-Ricci soliton}.
\end{definition}

In \cite{M2}, M\"uller showed that if $(g(t),u(t))$ is a gradient
harmonic-Ricci soliton with potential $f$, then
\begin{eqnarray*}
0&=&{\rm Ric}_{g(t)}-2du(t)\otimes du(t)+\nabla^{2}_{g(t)}f(t)
+cg(t), \\
0&=&\Delta_{g(t)}u(t)-\left\langle\nabla_{g(t)} u(t),
\nabla_{g(t)}f(t)\right\rangle_{g(t)}
\end{eqnarray*}
for some constant $c$.

\begin{corollary} \label{c6.6}There is no expanding breather on compact Riemannian manifolds other than expanding gradient harmonic-Ricci solitons.
\end{corollary}

\begin{proof} The proof is similar to that given in \cite{Li}. Suppose there is an expanding breather on a compact Riemannian manifold $M$, then by definition we have
\begin{equation*}
g(t_{2})=\alpha\Phi^{\ast}g(t_{1}), \ \ \ u(t_{2})=\Phi^{\ast}u(t_{1})
\end{equation*}
for some $t_{1}<t_{2}$, where $\Phi$ is a diffeomorphism and the
constant $\alpha>1$. Let $f_{+}(x)$ is a smooth function where
$\mathcal{W}_{+}(g(t_{2}),u(t_{2}),\tau(t_{2}),f(t_{2}))$ attains its
minimum. Then there exists a smooth function $f_{+}(x,t): M\times[t_{1},t_{2}]\to\mathbb{R}$ with initial value
$f_{+}(x,t_{2})=f_{+}(x)$ and satisfies the coupled system appeared
in \ref{r6.3}. Define a linear function
\begin{equation*}
\tau: [t_{1},t_{2}]\longrightarrow(0,+\infty), \ \ \ \tau(t_{2})=T+t_{2}
\end{equation*}
where $T$ is a constant. By the monotonicity formula, we have
\begin{eqnarray*}
\mu_{+}(g(t_{2}),u(t_{2}),\tau(t_{2}))&=&
\mathcal{W}_{+}(g(t_{2}),u(t_{2}),\tau(t_{2}),f_{+}(t_{2}))\\
&\geq&\mathcal{W}_{+}(g(t_{1}),u(t_{1}),\tau(t_{1}),f_{+}(t_{1})) \\
&\geq&\mu_{+}(g(t_{1}),u(t_{1}),\tau(t_{1})).
\end{eqnarray*}
Lemma \ref{l6.4} and the diffeomorphic invariant property of the
functionals shows
\begin{equation*}
\mu_{+}(g(t_{1}),u(t_{1}),\tau(t_{1}))\leq\alpha\mu_{+}(g(t_{1}),
u(t_{1}),\tau(t_{1}))
\end{equation*}
which yields
\begin{equation*}
\mu_{+}(g(t_{1}),u(t_{1}),\tau(t_{1}))\geq0
\end{equation*}
since $\alpha>1$.

If we impose an additional condition $\tau(t_{2})=\alpha\tau(t_{1})$ and $\tau(t_{1})=T+t_{1}$, we have
\begin{equation*}
\tau(t)=\frac{\alpha(t-t_{1})-(t-t_{2})}{\alpha-1}, \ \ \ T=\frac{t_{2}-\alpha t_{1}}{\alpha-1}.
\end{equation*}
Then
\begin{equation*}
\frac{\tau(t_{2})^{\frac{n}{2}}}{V_{g(t_{2})}}=\frac{\left[\frac{\alpha(t_{2}-t_{1})}{\alpha-1}\right]^{\frac{n}{2}}}
{\alpha^{\frac{n}{2}}V_{g(t_{1})}}=\frac{\tau(t_{1})^{\frac{n}{2}}}{V_{g(t_{1})}}.
\end{equation*}
The mean value theorem tells us that there exists a time $\overline{t}\in[t_{1},t_{2}]$ with
\begin{eqnarray*}
0&=&\frac{d}{dt}\Big|_{t=\overline{t}}{\rm log}\frac{\tau(t)^{\frac{n}{2}}}{V_{g(t)}}\\
&=&\frac{V_{g(\overline{t})}}{\tau(\overline{t})^{\frac{n}{2}}}\cdot\frac{\frac{n}{2}\tau(\overline{t})^{\frac{n}{2}-1}
V_{g(\overline{t})}-\tau(\overline{t})^{\frac{n}{2}}
\frac{d}{dt}\Big|_{t=\overline{t}}V_{g(t)}}{V^{2}_{g(\overline{t})}}\\
&=&\frac{n}{2\tau(\overline{t})}-\frac{1}{V_{g(\overline{t})}}\frac{\partial}{\partial t}\Big|_{t=\overline{t}}V_{g(\overline{t})}.
\end{eqnarray*}
From the evolution equation for the volume element $dV_{g(t)}$ we have
\begin{equation*}
\frac{d}{dt}V_{g(t)}=\int_{M}\frac{\partial}{\partial t}dV_{g(t)}
=\int_{M}\left(-S_{g(t),u(t)}-\Delta_{g(t)}f(t)\right)dV_{g(t)}
=-\int_{M}S_{g(t),u(t)}dV_{g(t)}.
\end{equation*}
Putting those together yields
\begin{equation*}
0=\frac{n}{2\tau(\overline{t})}+\frac{1}{V_{g(\overline{t})}}
\int_{M}S_{g(\overline{t}),u(\overline{t})}
dV_{g(\overline{t})}=\frac{1}{V_{g(\overline{t})}}\int_{M}
\left(S_{g(\overline{t}),u(\overline{t})}
+\frac{n}{2\tau(\overline{t})}\right)dV_{g(\overline{t})}.
\end{equation*}
If we set $\overline{f}={\rm log}V_{g(\overline{t})}$ then
\begin{equation*}
0=\mathcal{W}_{+}(g(\overline{t}),u(\overline{t}),\tau(\overline{t}),
\overline{f})
\geq\mu_{+}(g(\overline{t}),u(\overline{t}),\tau(\overline{t})).
\end{equation*}
By the monotonicity of $\mu_{+}$ we obtain
\begin{equation*}
0\leq\mu_{+}(g(t_{1}),u(t_{1}),\tau(t_{1}))\leq\mu_{+}(g(\overline{t}),u(\overline{t}),\tau(\overline{t}))
\leq0
\end{equation*}
Hence $\mu_{+}(g(t_{1}),u(t_{1}),\tau(t_{1}))=\mu_{+}(g(t_{2}),u(t_{2}),\tau(t_{2}))=0$ and $\mathcal{W}_{+}=0$ on the interval $[t_{1},t_{2}]$. This indicates that the first variation of $\mathcal{W}_{+}$ must vanish. So the expanding breather is a gradient soliton, i.e.,
\begin{equation*}
\mathcal{S}_{g(t),u(t)}+\nabla^{2}_{g(t)}f(t)+\frac{1}{2\tau(t)}g(t)=0.
\end{equation*}
Moreover, in this case $\Delta_{g(t)}u(t)=\left\langle du(t),df(t)\right\rangle_{g(t)}
$.
\end{proof}

As (\ref{6.7}), we define
\begin{equation}
\mu_{+,k}(g,u,\tau):=\inf\left\{\mathcal{W}_{+,k}(g,u,\tau,f)\Big| f\in C^{+\infty}(M), \ \ \ \int_{M}e^{-f}dV_{g}=1\right\}\label{6.8}
\end{equation}
As Lemma \ref{l6.4}, we still have
\begin{equation}
\mu_{+,k}(\alpha g,u,\alpha\tau)=\alpha\mu_{+,k}(g,u,\tau).\label{6.9}
\end{equation}

\begin{corollary} \label{c6.7}If $(g(t),u(t))$ is an expanding harmonic-Ricci breathers
on compact Riemannian manifolds, then $M$ is an Einstein manifold and $u(t)$
is constant.
\end{corollary}

\begin{proof} Using the same method in Corollary \ref{c6.6} and $\mu_{+,k}$, we can show
that the first variation of $\mathcal{W}_{+,k}$ must vanish. Hence, from (\ref{6.4}) one has
\begin{eqnarray*}
\mathcal{S}_{g(t),u(t)}+\nabla^{2}_{g(t)}f(t)+\frac{1}{2\tau(t)}g(t)&=&0,\\
\mathcal{S}_{g(t),u(t)}+\frac{1}{2\tau(t)}g(t)&=&0,\\
\Delta_{g(t)}u(t)&=&\left\langle du(t),df(t)\right\rangle_{g(t)},\\
\Delta_{g(t)}u(t)&=&0.
\end{eqnarray*}
The above four equations can be reduced to a coupled equation
\begin{equation*}
\mathcal{S}_{g(t),u(t)}+\frac{1}{2\tau(t)}g(t)=0=\Delta_{g(t)}u(t)
\end{equation*}
which indicates that $u(t)$ is a constant and ${\rm Ric}_{g(t)}=-\frac{1}{2\tau(t)}g(t)$.
\end{proof}

\section{Eigenvalues of the Laplacian under the harmonic-Ricci
flow}\label{section7}

In this section we consider the eigenvalues of the Laplacian $\Delta_{g(t)}$
under the harmonic-Ricci flow
\begin{eqnarray}
\frac{\partial}{\partial t}g(t)&=&-2{\rm Ric}_{g(t)}+4du(t)\otimes du(t),\label{7.1}\\
\frac{\partial}{\partial t}u(t)&=&\Delta_{g(t)}u(t).\label{7.2}
\end{eqnarray}
Suppose that $\lambda(t)$, which is a function of time $t$ only, is an eigenvalue of the Laplacian $\Delta_{g(t)}$ with an eigenfunction $f(t)=f(x,t)$, i.e.,
\begin{equation}
-\Delta_{g(t)}f(t)=\lambda(t)f(t).\label{7.3}
\end{equation}
Taking the derivative with respect to $t$, we get
\begin{equation*}
-\left(\frac{\partial}{\partial t}\Delta_{g(t)}\right)f(t)
-\Delta_{g(t)}\left(\frac{\partial}{\partial t}f(t)\right)=\left(\frac{d}{dt}\lambda(t)\right)f(t)+\lambda(t)\frac{\partial}{\partial t}f(t).
\end{equation*}
Integrating above equation with $f$ yields
\begin{eqnarray*}
&&-\int_{M}f(t)\left(\frac{\partial}{\partial t}\Delta_{g(t)}\right)f(t)\!\ dV_{g(t)}
-\int_{M}f(t)\Delta_{g(t)}\left(\frac{\partial}{\partial t}f(t)\right)\!\ dV_{g(t)}
\\
&=&\frac{d}{dt}\lambda(t)\cdot\int_{M}f(t)^{2}\!\ dV_{g(t)}+\lambda(t)\int_{M}f(t)\frac{\partial}{\partial t}f(t)\!\ dV_{g(t)}.
\end{eqnarray*}
Since
\begin{eqnarray*}
-\int_{M}f(t)\Delta\left(\frac{\partial}{\partial t}f(t)\right)dV_{g(t)}
&=&-\int_{M}\Delta_{g(t)}f(t)\cdot\frac{\partial}{\partial t}f(t)\!\ dV_{g(t)}\\
&=&\lambda(t)\int_{M}f(t)\frac{\partial}{\partial t}f(t)\!\ dV_{g(t)},
\end{eqnarray*}
it follows that
\begin{equation}
\frac{d}{dt}\lambda(t)\cdot\int_{M}f(t)^{2}\!\ dV_{g(t)}
=-\int_{M}f(t)\left(\frac{\partial}{\partial t}\Delta_{g(t)}\right)
f(t)\!\ dV_{g(t)}.\label{7.4}
\end{equation}
If we set $v_{ij}=-2R_{ij}+4\partial_{i}u\partial_{j}u$, then
\begin{equation*}
\frac{\partial}{\partial t}\Gamma^{k}_{ij}=\frac{1}{2}g^{k\ell}
\left(\partial_{i}v_{\ell j}+\partial_{j}v_{il}-\partial_{\ell}v_{ij}\right).
\end{equation*}
We temporarily omit all subscripts $t$. Multiplying with $g^{ij}$ on both sides, we obtain
\begin{eqnarray*}
g^{ij}\frac{\partial}{\partial t}\Gamma^{k}_{ij}&=&\frac{1}{2}g^{kl}\left(2\nabla^{i}v_{li}-\nabla_{l}\left(g^{ij}v_{ij}\right)
\right) \ \ = \ \ g^{kl}\nabla^{i}v_{il}+\nabla^{k}S\\
&=&g^{kl}\nabla^{i}\left(-2R_{il}+4\nabla_{i}u\nabla_{l}u\right)+\nabla^{k}\left(
R-2|\nabla u|^{2}\right)\\
&=&-\nabla^{k}R+4\Delta u\cdot\nabla^{k}u+4\nabla_{i}u\cdot\nabla^{i}\nabla^{k}u
+\nabla^{k}R-4\nabla^{k}\nabla^{i}u\cdot\nabla_{i}u\\
&=&4\Delta u\cdot\nabla^{k}u.
\end{eqnarray*}
Therefore,
\begin{eqnarray*}
\frac{\partial}{\partial t}\left(\Delta f\right)&=&
\frac{\partial}{\partial t}\left(g^{ij}\nabla_{i}\nabla_{j}f\right)\\
&=&\left(\frac{\partial}{\partial t}g^{ij}\right)\nabla_{i}\nabla_{j}f
+g^{ij}\left[\partial_{i}\partial_{j}\frac{\partial f}{\partial t}-\left(\frac{\partial}{\partial t}\Gamma^{k}_{ij}\right)\partial_{k}f
-\Gamma^{k}_{ij}\partial_{k}\frac{\partial f}{\partial t}\right]\\
&=&\left(\frac{\partial}{\partial t}g^{ij}\right)\nabla_{i}\nabla_{j}f
+\Delta_{g(t)}\left(\frac{\partial}{\partial t}f\right)-g^{ij}\left(\frac{\partial}{\partial t}\Gamma^{k}_{ij}\right)\nabla_{k}f\\
&=&\left(2R_{ij}-4\nabla_{i}u\nabla_{j}u\right)\nabla^{i}\nabla^{j}f-4\Delta u\cdot\nabla^{k}u\nabla_{k}f+\Delta_{g(t)}\left(\frac{\partial}{\partial t}f\right).
\end{eqnarray*}
Plugging it into (\ref{7.4}) we derive
\begin{eqnarray*}
\frac{d}{dt}\lambda(t)\cdot\int_{M}f(t)^{2}\!\ dV_{g(t)}&=&-2\int_{M}R_{ij}\nabla^{i}\nabla^{j}f\!\ dV+4\int_{M}f\nabla^{i}u\nabla^{j}u
\nabla_{i}\nabla_{j}f\!\ dV\\
&&+ \ 4\int_{M}f\Delta u\cdot\nabla^{k}u\nabla_{k}f\!\ dV.
\end{eqnarray*}
The first term can be rewritten as
\begin{eqnarray*}
-2\int_{M}fR_{ij}\nabla^{i}\nabla^{j}f\!\ dV&=&\int_{M}\nabla^{i}\left(2fR_{ij}\right)\nabla^{j}fdV\\
&=&2\int_{M}\left(\nabla^{i}f\cdot R_{ij}+f\cdot\nabla^{i}R_{ij}\right)\nabla^{j}f\!\ dV\\
&=&2\int_{M}R_{ij}\nabla^{i}f\nabla^{j}f\!\ dV+\int_{M}f\nabla_{j}R\nabla^{j}f\!\ dV\\
&=&2\int_{M}R_{ij}\nabla^{i}f\nabla^{j}f\!\ dV-\int_{M}R\nabla_{j}\left(f\nabla^{j}f\right)\!\ dV\\
&=&\lambda\int_{R}f^{2}\!\ dV-\int_{M}R|\nabla f|^{2}\!\ dV+2\int_{M}R_{ij}\nabla^{i}f\nabla^{j}f\!\ dV.
\end{eqnarray*}
Hence
\begin{eqnarray*}
\left(\frac{d}{dt}\lambda(t)\right)\int_{M}f(t)^{2}
\!\ dV_{g(t)}&=&
\lambda(t)\int_{M}R_{g(t)}f(t)^{2}dV_{g(t)}+2\int_{M}R_{ij}\nabla^{i}f\nabla^{j}f\!\ dV\\
&&- \ \int_{M}R_{g(t)}\left|\nabla_{g(t)} f(t)\right|^{2}_{g(t)}dV_{g(t)}
\\
&&+ \ 4\int_{M}f\left(\nabla^{i}u\nabla^{j}u\nabla_{i}
\nabla_{j}f+\Delta u\nabla^{k}u\nabla_{k}f\right)\!\ dV.
\end{eqnarray*}
On the other hand,
\begin{eqnarray*}
&&\int_{M}f\nabla^{i}u\nabla^{j}u\nabla_{i}\nabla_{j}f\!\ dV \ \ = \ \ -\int_{M}
\nabla_{i}\left(f\nabla^{i}u\nabla^{j}u\right)\nabla_{j}f\!\ dV\\
&=&-\int_{M}\left(\nabla_{i}f\nabla^{i}u\nabla^{j}u+f\Delta u\nabla^{j}u
+f\nabla^{i}u\nabla_{i}\nabla^{j}u\right)\nabla_{j}f\!\ dV\\
&=&-\int_{M}f\Delta u\langle \nabla u,\nabla f\rangle dV
-\int_{M}\nabla^{i}u\nabla^{j}u\nabla_{i}f\nabla_{j}f\!\ dV\\
&&- \ \int_{M}f\nabla^{i}u\nabla^{j}f
\nabla_{i}\nabla_{j}u\!\ dV
\end{eqnarray*}
and therefore
\begin{eqnarray*}
\frac{d}{dt}\lambda(t)\int_{M}f(t)^{2}dV_{g(t)}
&=&\lambda(t)\int_{M}R_{g(t)}f(t)^{2}dV_{g(t)}
-4\int_{M}f\nabla^{i}u\nabla^{j}f\nabla_{i}\nabla_{j}u\!\ dV\\
&&+ \ 2\int_{M}S_{ij}\nabla^{i}f\nabla_{j}f\!\ dV-\int_{M}R_{g(t)}\left|\nabla_{g(t)} f(t)\right|^{2}_{g(t)}dV_{g(t)}.
\end{eqnarray*}
The last term in above can be simplified as follows:
\begin{eqnarray*}
&&-\int_{M}f\nabla^{i}u\nabla^{j}f\nabla_{i}\nabla_{j}u\!\ dV \ \ = \ \ \int_{M}
\nabla^{j}\left(f\nabla_{i}u\nabla_{j}f\right)\nabla^{i}u\!\ dV\\
&=&\int_{M}\left(\nabla^{j}f\nabla_{i}u\nabla_{j}f+f\nabla^{j}\nabla_{i}u\nabla_{j}f
+f\nabla_{i}u\Delta f\right)\nabla^{i}u\!\ dV\\
&=&\int_{M}|\nabla u|^{2}|\nabla f|^{2}dV+\int_{M}f\Delta f|\nabla u|^{2}dV+\int_{M}f\nabla^{i}u\nabla^{j}f\nabla_{i}\nabla_{j}u\!\ dV
\end{eqnarray*}
consequently,
\begin{equation*}
-2\int_{M}f\nabla^{i}u\nabla^{j}f\nabla_{i}\nabla_{j}u\!\ dV
=\int_{M}|\nabla u|^{2}|\nabla f|^{2}dV-\lambda\int_{M}f^{2}|\nabla u|^{2}dV.
\end{equation*}

Therefore we derive the following

\begin{theorem}\label{t7.1}If $(g(t),u(t))$ is a solution of the harmonic-Ricci flow on a compact Riemannian manifold $M$ and $\lambda(t)$ denotes the eigenvalue of the
Laplacian $\Delta_{g(t)}$, then
\begin{eqnarray}
\frac{d}{dt}\lambda(t)\cdot\int_{M}f(t)^{2}dV_{g(t)}&=&\lambda(t)
\int_{M}S_{g(t),u(t)}f(t)^{2}dV_{g(t)}\nonumber\\
&&-
\int_{M}S_{g(t),u(t)}\left|\nabla_{g(t)} f(t)\right|^{2}_{g(t)}dV_{g(t)}\label{7.5}\\
&&+2\int_{M}\left\langle\mathcal{S}_{g(t),u(t)},df(t)\otimes df(t)\right\rangle dV_{g(t)}.\nonumber
\end{eqnarray}
\end{theorem}

We set
\begin{equation}
S_{\min}(0):=\min_{x\in M}S(x,0).\label{7.6}
\end{equation}

\begin{theorem} \label{t7.2}Let $(g(t),u(t))_{t\in[0,T]}$ be a solution of the harmonic-Ricci flow on a compact Riemannian manifold $M$ and $\lambda(t)$ denote the eigenvalue of the Laplacian $\Delta_{g(t)}$. Suppose that $\mathcal{S}_{g(t),u(t)}-\alpha S_{g(t),u(t)}g(t)\geq0$ along the harmonic-Ricci flow for some $\alpha\geq\frac{1}{2}$.

\begin{itemize}

\item[(1)] If $S_{\min}(0)\geq0$, then $\lambda(t)$ is nondecreasing along the harmonic-Ricci flow for any $t\in[0,T]$.

\item[(2)] If $S_{\min}(0)>0$, then the quantity
\begin{equation*}
\left(1-\frac{2}{n}S_{\min}(0)t\right)^{n\alpha}\lambda(t)
\end{equation*}
is nondecreasing along the harmonic-Ricci flow for
$T\leq\frac{n}{2S_{\min}(0)}$.

\item[(3)] If $S_{\min}(0)<0$, then the quantity
\begin{equation*}
\left(1-\frac{2}{n}S_{\min}(0)t\right)^{n\alpha}\lambda(t)
\end{equation*}
is nondecreasing along the harmonic-Ricci flow for any $t\in[0,T]$.

\end{itemize}
\end{theorem}

\begin{proof} By Theorem \ref{t7.1}, we have
\begin{eqnarray*}
\frac{d}{dt}\lambda(t)&\geq&\left(\frac{\int_{M}S_{g(t),u(t)}f(t)^{2}
dV_{g(t)}}{\int_{M}f(t)^{2}dV_{g(t)}}\right)\lambda(t)\\
&&+(2\alpha-1)\left(\frac{\int_{M}S_{g(t),u(t)}
|\nabla_{g(t)} f(t)|^{2}_{g(t)}}{\int_{M}f(t)^{2}dV_{g(t)}}\right).
\end{eqnarray*}
By definition we have $-f(t)\Delta_{g(t)}=\lambda(t)f(t)$. Taking the integration on both sides yields that $\lambda(t)\geq0$. Since
\begin{equation*}
\frac{\partial}{\partial t}S_{g(t),u(t)}=\Delta_{g(t)}S_{g(t),u(t)}
+2\left|\mathcal{S}_{g(t),u(t)}\right|^{2}_{g(t)}+4\left|\Delta_{g(t)}u(t)\right|^{2}_{g(t)}
\end{equation*}
and $\left|\mathcal{S}_{g(t),u(t)}\right|^{2}\geq\frac{1}{n}S_{g(t),u(t)}^{2}$, it
follows that
\begin{equation*}
\frac{\partial}{\partial t}S_{g(t),u(t)}\geq\Delta_{g(t)}S_{g(t),u(t)}+\frac{2}{n}S^{2}_{g(t),u(t)}.
\end{equation*}
The corresponding ODE
\begin{equation*}
\frac{d}{dt}a(t)=\frac{2}{n}a(t)^{2}, \ \ \ a(t)=S_{\min}(0)
\end{equation*}
has the solution
\begin{equation*}
a(t)=\frac{S_{\min}(0)}{1-\frac{2}{n}S_{\min}(0)t}.
\end{equation*}
Then the maximum principle implies $S_{g(t),u(t)}\geq a(t)$ and hence, using the assumption that $2\alpha-1\geq0$,
\begin{equation*}
\frac{d}{dt}\lambda(t)\geq a(t)\lambda(t)+(2\alpha-1)a(t)\frac{\int_{M}|\nabla_{g(t)} f(t)|^{2}_{g(t)}dV_{g(t)}}{\int_{M}f(t)^{2}dV_{g(t)}}.
\end{equation*}
By integration by parts, we note that
\begin{equation*}
\int_{M}\left|\nabla f\right|^{2}dV=-\int_{M}f\cdot\Delta fdV=\lambda\int_{M}f^{2}\!\ dV
\end{equation*}
which indicates
\begin{equation*}
\frac{d}{dt}\lambda(t)\geq a(t)\lambda(t)+(2\alpha-1)a(t)\lambda=2\alpha a(t)\lambda(t)
\end{equation*}
and
\begin{equation*}
\frac{d}{dt}\left(\lambda(t)\cdot \exp\left(-2\alpha\int^{t}_{0}a(\tau)d\tau\right)\right)\geq0.
\end{equation*}
Plugging the expression into above yields the desired result. If $S_{\min}(0)\geq0$, by the nonnegativity of $\mathcal{S}_{g(t)}$ preserved along the harmonic-Ricci flow, we conclude that $\frac{d}{dt}\lambda(t)\geq0$.
\end{proof}

\begin{corollary} \label{c7.3}Let $(g(t),u(t))_{t\in[0,T]}$ be a solution of the harmonic-Ricci flow on a compact Riemannian surface $\Sigma$ and $\lambda(t)$ denote the eigenvalue of the Laplacian $\Delta_{g(t)}$.

\begin{itemize}

\item[(1)] Suppose that ${\rm Ric}_{g(t)}\leq\epsilon du(t)\otimes du(t)$ where
\begin{equation}
\epsilon\leq4\frac{1-\alpha}{1-2\alpha}, \ \ \ \alpha>\frac{1}{2}.\label{10.7}
\end{equation}

\begin{itemize}

\item[(1-1)] If $S_{\min}(0)\geq0$, then $\lambda(t)$ is nondecreasing along the harmonic-Ricci flow for any $t\in[0,T]$.

\item[(1-2)] If $S_{\min}(0)>0$, then the quantity
\begin{equation*}
\left(1-S_{\min}(0)t\right)^{2\alpha}\lambda(t)
\end{equation*}
is nondecreasing along the harmonic-Ricci flow for
$T\leq\frac{1}{S_{\min}(0)}$.

\item[(1-3)] If $S_{\min}(0)<0$, then the quantity
\begin{equation*}
\left(1-S_{\min}(0)t\right)^{2\alpha}\lambda(t)
\end{equation*}
is nondecreasing along the harmonic-Ricci flow for any $t\in[0,T]$.

\end{itemize}

\item[(2)] Suppose that
\begin{equation}
\left|\nabla_{g(t)} u(t)\right|^{2}_{g(t)}g(t)\geq2du(t)\otimes du(t).\label{10.8}
\end{equation}

\begin{itemize}

\item[(2-1)] If $S_{\min}(0)\geq0$, then $\lambda(t)$ is nondecreasing along the harmonic-Ricci flow for any $t\in[0,T]$.

\item[(2-2)] If $S_{\min}(0)>0$, then the quantity
\begin{equation*}
\left(1-S_{\min}(0)t\right)\lambda(t)
\end{equation*}
is nondecreasing along the harmonic-Ricci flow for
$T\leq\frac{1}{S_{\min}(0)}$.

\item[(2-3)] If $S_{\min}(0)<0$, then the quantity
\begin{equation*}
\left(1-S_{\min}(0)t\right)\lambda(t)
\end{equation*}
is nondecreasing along the harmonic-Ricci flow for any $t\in[0,T]$.

\end{itemize}
\end{itemize}
\end{corollary}

\begin{proof} As above, we always omit subscripts $t$.
In the surface case, we have $R_{ij}=\frac{R}{2}g_{ij}$. Then
\begin{eqnarray*}
T_{ij}&:=&S_{ij}-\alpha Sg_{ij} \ \ = \ \ \frac{R}{2}g_{ij}-2\nabla_{i}u\nabla_{j}u-\alpha\left(R-2\left|\nabla u\right|^{2}\right)g_{ij}\\
&=&\left(\frac{1}{2}-\alpha\right)Rg_{ij}-2\nabla_{i}u\nabla_{j}u+2\alpha\left|\nabla u\right|^{2}g_{ij}.
\end{eqnarray*}
For any vector $V=(V^{i})$, we calculate
\begin{eqnarray*}
T_{ij}V^{i}V^{j}&=&\left(\frac{1}{2}-\alpha\right)R|V|^{2}-2\left(\nabla_{i}u V^{i}\right)^{2}+2\alpha|\nabla u|^{2}|V|^{2}\\
&\geq&\left(\frac{1}{2}-\alpha\right)R|V|^{2}-2|\nabla u|^{2}|V|^{2}+2\alpha|\nabla u|^{2}|V|^{2}.
\end{eqnarray*}
If $R_{ij}\leq\epsilon \nabla_{i}u\nabla_{j}u$, then $T_{ij}V^{i}V^{j}=\left[\left(\frac{1}{2}-\alpha\right)\epsilon-2+2\alpha\right]|\nabla u|^{2}|V|^{2}\geq0$.

For the second case, we note that
\begin{eqnarray*}
T_{ij}V^{i}V^{j}&=&R_{ij}V^{i}V^{j}-2\nabla_{i}u V^{i}\nabla_{j}u V^{j}-\frac{R}{2}|V|^{2}+|\nabla u|^{2}|V|^{2}\\
&\geq&R_{ij}V^{i}V^{j}-|\nabla u|^{2}|V|^{2}-\frac{R}{2}|V|^{2}+|\nabla u|^{2}|V|^{2} \ \ = \ \ 0.
\end{eqnarray*}
Hence, the corresponding results follow by Theorem \ref{t7.2}.
\end{proof}

When we consider the Ricci flow, we have the following two results derived from Corollary \ref{c7.3}.

\begin{corollary} \label{c7.4}Let $(g(t))_{t\in[0,T]}$ be a solution of the Ricci flow on a compact Riemannian surface $\Sigma$ and $\lambda(t)$ denote the eigenvalue of the Laplacian $\Delta_{g(t)}$.

\begin{itemize}

\item[(1)] If $R_{\min}(0)\geq0$, then $\lambda(t)$ is nondecreasing along the Ricci flow for any $t\in[0,T]$.

\item[(2)] If $R_{\min}(0)>0$, then the quantity $(1-R_{\min}(0)t)\lambda(t)$ is nondecreasing along the Ricci flow for $T\leq\frac{1}{R_{\min}(0)}$.

\item[(3)] If $R_{\min}(0)<0$, then the quantity $(1-R_{\min}(0)t)\lambda(t)$ is nondecreasing along the Ricci flow for any $t\in[0,T]$.

\end{itemize}
\end{corollary}

\begin{remark} \label{r7.5}Let $(g(t))_{t\in[0,T]}$ be a solution of the Ricci flow on a compact Riemannian surface $\Sigma$ with nonnegative scalar curvature and $\lambda(t)$ denote the eigenvalue of the Laplacian $\Delta_{g(t)}$. Then $\lambda(t)$ is nondecreasing along the Ricci flow for $t\in[0,T]$.
\end{remark}

\section{Eigenvalues of the Laplacian-type under the harmonic-Ricci flow}\label{section8}

Recall that
\begin{equation}
\mu(g,u)=\mu_{1}(g,u)=\inf\left\{\mathcal{F}(g,u,f)\Big|\int_{M}e^{-f}dV_{g}
=1\right\}.\label{8.1}
\end{equation}
We showed that $\mu(g,u)$ is the smallest eigenvalue of the operator $-4\Delta_{g}+R_{g}-2|\nabla_{g} u|^{2}_{g}$. Inspired by \cite{C1,C2}, we define a Laplacian-type operators associated with quantities $g,u,c$:
\begin{eqnarray}
\Delta_{g,u,c}&:=&-\Delta_{g}+c\left(R_{g}-2\left|\nabla_{g} u\right|^{2}_{g}\right), \label{8.2}\\
\Delta_{g,u}&:=&\Delta_{g,u,\frac{1}{2}} \ \ = \ \ -\Delta_{g}+\frac{1}{2}\left(R_{g}-2\left|\nabla_{g} u\right|^{2}_{g}\right).\label{8.3}
\end{eqnarray}
Then $\mu(g,u)$ is the smallest eigenvalue of the operator $4\Delta_{g,u,\frac{1}{4}}$.

To the operator $\Delta_{g,u}$ we associate a functional
\begin{equation}
C^{\infty}(M)\longrightarrow\mathbb{R}, \ \ \ f\longmapsto\lambda_{g,u}(f):=\int_{M}f\Delta_{g,u}f\!\ dV_{g}.\label{8.4}
\end{equation}
When $f$ is an eigenfunction of the operator $\Delta_{g,u}$ with the eigenvalue $\lambda$, i.e., $\Delta_{g,u}f=\lambda f$ and normalized by $\int_{M}f^{2}dV_{g}=1$, we obtain $\lambda_{g,u}(f)
=\lambda$. Next lemma will deal with the evolution equation for $\lambda(f(t))$ where $f(t)$ is an eigenfunction of $\Delta_{g(t),u(t)}$ and the couple $(g(t),u(t))$ satisfies the harmonic-Ricci flow. Set
\begin{equation}
v_{ij}:=-2S_{ij}=-2R_{ij}+4\partial_{i}u\cdot\partial_{j}u, \ \ \ v:=g^{ij}v_{ij}.\label{4.9}
\end{equation}
The obtained symmetric tensor field is denoted by $\mathcal{V}_{g(t),u(t)}=(v_{ij})$.

\begin{lemma} \label{l8.1}Suppose that $(g(t),u(t))$ is a solution of the harmonic-Ricci flow on a compact Riemannian manifold $M$ and $f(t)$ is an eigenfunction of $\Delta_{g(t),u(t)}$, i.e., $\Delta_{g(t),u(t)}f(t)=\lambda(t)f(t)$(where $\lambda(t)$
is only a function of time $t$ only), with the normalized condition $\int_{M}f(t)^{2}dV_{g(t)}=1$. Then we have
\begin{eqnarray}
&&\frac{d}{dt}\lambda_{g(t),u(t)}(f(t)\nonumber\\
&=&\int_{M}f(t) \left(\nabla^{i}v_{ik}-\frac{1}{2}\nabla_{k}v\right)
\nabla^{k}f(t)dV_{g(t)}-\int_{M}f^{2}(t)\frac{\partial}{\partial t}\left|\nabla_{g(t)} u(t)\right|^{2}_{g(t)}dV_{g(t)}
\nonumber\\
&&+ \ \int_{M}\left[\left\langle\mathcal{V}_{g(t),u(t)},
\nabla^{2}_{g(t)}f(t)\right\rangle_{g(t)}+\frac{1}{2}\left(\frac{\partial}{\partial t}R_{g(t)}\right)f(t)\right]f(t)dV_{g(t)}.\label{8.6}
\end{eqnarray}
\end{lemma}

Before proving the lemma, we recall a formula that is an immediate consequence
of the evolution equation:
\begin{eqnarray}
\frac{\partial}{\partial t}\left(\Delta_{g(t)}f\right)&=&-g^{ip}g^{jq}v_{pq}
\nabla_{i}\nabla_{j}f-g^{ij}g^{k\ell}\nabla_{i}v_{j\ell}
\nabla_{k}f\label{8.7}\\
&&+\frac{1}{2}\left\langle\nabla_{g(t)} v_{g(t)},\nabla_{g(t)} f(t)\right\rangle_{g(t)}\nonumber
\end{eqnarray}
where the metric $g(t)$ evolves by $\frac{\partial}{\partial t}g_{ij}=v_{ij}$.

\begin{proof} Using (\ref{8.7}) and integration by parts, we get
\begin{eqnarray*}
&&\frac{d}{dt}\lambda_{g(t),u(t)}(f(t))\\
&=&\frac{\partial}{\partial t}\int_{M}\left[-\Delta_{g(t)}f(t)+\left(\frac{R_{g(t)}}{2}
-\left|\nabla_{g(t)} u(t)\right|^{2}_{g(t)}\right)f(t)\right]f(t)dV_{g(t)}\\
&=&\int_{M}\left[g^{ip}g^{jq}v_{pq}\nabla_{i}\nabla_{j}f
+g^{ij}g^{kl}\nabla_{i}v_{jl}\nabla_{k}f
-\frac{1}{2}\left\langle\nabla_{g(t)} v_{g(t)},\nabla_{g(t)} f(t)\right\rangle_{g(t)}\right]\\
&&f(t)dV_{g(t)}+\int_{M}\left[-\Delta_{g(t)}\left(\frac{\partial}{\partial t}f(t)\right)+\left(\frac{R_{g(t)}}{2}-\left|\nabla_{g(t)} u(t)\right|^{2}_{g(t)}\right)\frac{\partial}{\partial t}f(t)\right.\\
&& \ +\left.\left(\frac{\partial}{\partial t}\left(\frac{1}{2}R_{g(t)}\right)-\frac{\partial}{\partial t}\left(\left|\nabla_{g(t)} u(t)\right|^{2}_{g(t)}\right)\right)f(t)\right]f(t)dV_{g(t)}\\
&&+ \ \int_{M}\left[-\Delta_{g(t)}f(t)+\left(\frac{R_{g(t)}}{2}-\left|\nabla_{g(t)} u\right|^{2}_{g(t)}\right)f(t)\right]\frac{\partial}{\partial t}\left(f(t)dV_{g(t)}\right)\\
&=&\int_{M}\left(g^{ip}g^{jq}v_{pq}\nabla_{i}\nabla_{j}f+\frac{1}{2}\left(\frac{\partial}{\partial t}R_{g(t)}\right)f(t)\right)f(t)dV_{g(t)}\\
&&+ \ \int_{M}\left(g^{ij}g^{kl}\nabla_{i}v_{jl}\nabla_{k}f-\frac{1}{2}g^{kl}\nabla_{l}v\nabla_{k}f\right)
f(t)dV_{g(t)}\\
&&+ \ \int_{M}\Delta_{g(t),u(t)}f(t)\left(\frac{\partial}{\partial t}f(t)dV_{g(t)}+\frac{\partial}{\partial t}\left(f(t)dV_{g(t)}\right)\right)\\
&&- \ \int_{M}\frac{\partial}{\partial t}\left(
\left|\nabla_{g(t)} u(t)\right|^{2}_{g(t)}\right)f(t)^{2}dV_{g(t)}.
\end{eqnarray*}
Since $f(t)$ is an eigenfunction of $\Delta_{g(t),u(t)}$, it follows that
\begin{equation*}
\int_{M}\Delta_{g(t),u(t)}f(t)\left(\frac{\partial}{\partial t}f(t)dV_{g(t)}+\frac{\partial}{\partial t}\left(f(t)dV_{g(t)}\right)\right)
=\lambda(t)\frac{\partial}{\partial t}\int_{M}f(t)^{2}
dV_{g(t)}= 0
\end{equation*}
by the normalized condition. Thus we complete the proof.
\end{proof}

Using (\ref{3.6}), we find that the first term in the right hand side of (\ref{8.6}) can be written as
\begin{eqnarray*}
&&\int_{M}\left[v_{ij}\nabla^{i}\nabla^{j}f+\frac{1}{2}\left(\frac{\partial}{\partial t}R_{g(t)}\right)f(t)\right]f(t)dV_{g(t)}\\
&=&\int_{M}\left[-2f(t)\left\langle{\rm Ric}_{g(t)},\nabla^{2}_{g(t)}f(t)\right\rangle_{g(t)}\right.\\
&& \ +\left.4f(t)\left\langle
du(t)\otimes du(t),\nabla^{2}_{g(t)}f(t)\right\rangle_{g(t)}\right]dV_{g(t)}\\
&&+ \ \int_{M}\left[\left(\frac{1}{2}\Delta_{g(t)}R_{g(t)}+\left|{\rm Ric}_{g(t)}\right|^{2}_{g(t)}\right)f(t)^{2}+
2f(t)^{2}\left|\Delta_{g(t)}u(t)\right|^{2}_{g(t)}
\right.\\
&& \ -\left.2f(t)^{2}\left|\nabla^{2}_{g(t)} u(t)\right|^{2}_{g(t)}
-4f(t)^{2}\left\langle{\rm Ric}_{g(t)},du(t)\otimes du(t)\right\rangle_{g(t)} \right]dV_{g(t)}\\
&=&\int_{M}\left[-2f(t)\left\langle{\rm Ric}_{g(t)},\nabla^{2}_{g(t)}f(t)\right\rangle_{g(t)}\right.\\
&& \ +\left.\left(\frac{1}{2}\Delta_{g(t)}R_{g(t)}+\left|{\rm Ric}_{g(t)}\right|^{2}_{g(t)}\right)f(t)^{2}\right]dV_{g(t)}\\
&&+ \ \int_{M}\left[4f(t)\left\langle du\otimes du,\nabla^{2}_{g(t)}f(t)\right\rangle_{g(t)}-4f^{2}\left\langle du(t)\otimes du(t),{\rm Ric}_{g(t)}\right\rangle_{g(t)}\right.\\
&& \ +\left.2f(t)^{2}\left|\Delta_{g(t)}u(t)\right|^{2}_{g(t)}
-2f(t)^{2}\left|\nabla^{2}_{g(t)}u(t)\right|^{2}_{g(t)}\right]
dV_{g(t)}
\end{eqnarray*}
For the second term in (\ref{8.6}) one has, using the contracted
Bianchi identities,
\begin{eqnarray*}
&&\int_{M}\left(g^{ij}\nabla_{i}v_{jk}-\frac{1}{2}\nabla_{k}v\right)\nabla^{k}f\cdot f(t)\!\ dV_{g(t)}\\
&=&\int_{M}\bigg[g^{ij}\nabla_{i}\left(-2R_{jk}
+4\partial_{j}u\partial_{k}u\right)\\
&&- \ \frac{1}{2}
\nabla_{k}\left(-2R_{g(t)}+4\left|\nabla_{g(t)} u(t)\right|^{2}_{g(t)}\right)\bigg]\nabla^{k}f\cdot f(t)dV_{g(t)}\\
&=&\int_{M}4f(t)\Delta_{g(t)}u(t)\left\langle\nabla_{g(t)} u(t),
\nabla_{g(t)} f(t)\right\rangle_{g(t)}dV_{g(t)}\\
&&+\int_{M}\left(4g^{ij}\nabla_{j}u\cdot\nabla_{i}\nabla_{k}u
-2\nabla_{k}
\left|\nabla_{g(t)} u(t)\right|^{2}_{g(t)}\right)
\nabla^{k}f\cdot f(t)\!\ dV_{g(t)}\\
&=&\int_{M}4f(t)\Delta_{g(t)}u(t)\left\langle\nabla_{g(t)} u(t),\nabla_{g(t)} f(t)\right\rangle_{g(t)}dV_{g(t)}
\end{eqnarray*}
where in the last step we use the identity $\nabla_{k}|\nabla u|^{2}
=2g^{pq}\nabla_{k}\nabla_{p}u\cdot\nabla_{q}u$. Therefore
\begin{eqnarray}
&&\frac{d}{dt}\lambda_{g(t),u(t)}(f(t)) \ \ = \ \ \int_{M}\left[-2f(t)\left\langle{\rm Ric}_{g(t)},\nabla^{2}_{g(t)}f(t)\right\rangle_{g(t)}\right.\nonumber\\
&&\left.+\left(\frac{1}{2}\Delta_{g(t)}R_{g(t)}+\left|{\rm Ric}_{g(t)}\right|^{2}_{g(t)}\right)f(t)^{2}\right]dV_{g(t)}\nonumber\\
&&+ \ 4f(t)\int_{M}\left[\left\langle du(t)\otimes du(t),\nabla^{2}_{g(t)}f(t)\right\rangle_{g(t)}\right.-f(t)\left\langle du(t)\otimes du(t),{\rm Ric}_{g(t)}\right\rangle_{g(t)}\label{8.8}\\
&&+2f(t)^{2}\left|\Delta_{g(u)}u(t)\right|^{2}_{g(t)}
-2f(t)^{2}\left|\nabla^{2}_{g(t)}u(t)\right|^{2}_{g(t)}\nonumber\\
&& \ +\left.4f(t)\Delta_{g(t)}u(t)\left\langle\nabla_{g(t)} u(t),\nabla_{g(t)} f(t)\right\rangle_{g(t)}\right]
dV_{g(t)}\nonumber\\
&&- \ \int_{M}\left(\Delta_{g(t)}\left|\nabla_{g(t)}u(t)\right|^{2}_{g(t)}
-2\left|\nabla^{2}_{g(t)}u(t)\right|^{2}_{g(t)}-4\left|\nabla_{g(t)} u(t)\right|^{4}_{g(t)}\right)f(t)^{2}dV_{g(t)}.\nonumber
\end{eqnarray}

The above evolution equation can be simplified as

\begin{theorem} \label{t8.2}Suppose $(g(t),u(t))$ is a solution of the harmonic-Ricci flow on a compact Riemannian manifold $M$ and $f(t)$ is an eigenfunction of $\Delta_{g(t),u(t)}$, i.e., $\Delta_{g(t),u(t)}f(t)=\lambda(t)f(t)$(where $\lambda(t)$
is only a function of time $t$ only), with the normalized condition $\int_{M}f(t)^{2}dV_{g(t)}=1$. Then we have
\begin{eqnarray}
\frac{d}{d t}\lambda_{g(t),u(t)}(f(t))
&=&\int_{M}2\left\langle\mathcal{S}_{g(t)},df(t)\otimes
df(t)\right\rangle_{g(t)}dV_{g(t)}\nonumber\\
&&+\int_{M}f(t)^{2}\left[\left|\mathcal{S}_{g(t)}\right|^{2}_{g(t)}
\label{8.9}+2\left|\Delta_{g(t)}u(t)\right|^{2}_{g(t)}
\right]dV_{g(t)}.
\end{eqnarray}
\end{theorem}

\begin{proof} Calculate
\begin{eqnarray*}
&&\int_{M}4f(t)\Delta_{g(t)}u(t)\left\langle\nabla_{g(t)} u(t),\nabla_{g(t)} f(t)\right\rangle_{g(t)}dV_{g(t)}\\
&=&-4\int_{M}\nabla_{i}u\left[\nabla^{i}f\cdot\langle\nabla u,\nabla f
\rangle+f\left(\nabla^{i}\langle\nabla u,\nabla f\rangle\right)\right]dV\\
&=&-4\int_{M}\left|\langle\nabla u,\nabla f\rangle\right|^{2}dV_{g}-4\int_{M}f\nabla_{i}u\left(\langle\nabla^{i}\nabla u,\nabla f\rangle+
\langle\nabla u,\nabla^{i}\nabla f\right)dV.
\end{eqnarray*}
By the same method, we have
\begin{eqnarray*}
&&\int_{M}-\Delta_{g(t)}\left|\nabla_{g(t)} u(t)\right|^{2}_{g(t)}f(t)^{2}dV_{g(t)} \ \ = \ \ -\int_{M}|\nabla u|^{2}\left(2f\Delta f+2|\nabla f|^{2}\right)dV\\
&=&-2\int_{M}|\nabla f|^{2}|\nabla u|^{2}dV
-2\int_{M}f\Delta f|\nabla u|^{2}dV.
\end{eqnarray*}
However,
\begin{eqnarray*}
&&\int_{M}f\Delta f|\nabla u|^{2}dV \ \ = \ \ \int_{M}-\nabla_{i}f\cdot\nabla^{i}\left(f|\nabla u|^{2}\right)dV\\
&=&-\int_{M}\nabla_{i}f\left(\nabla^{i}f|\nabla u|^{2}+f\nabla^{i}|\nabla u|^{2}\right)dV\\
&=&-\int_{M}|\nabla u|^{2}|\nabla f|^{2}dV-\int_{M}f\nabla_{i}f\cdot\nabla^{i}|\nabla u|^{2}dV.
\end{eqnarray*}
Plugging it into above yields
\begin{eqnarray*}
&&\int_{M}-\Delta_{g(t)}\left|\nabla_{g(t)} u(t)\right|^{2}_{g(t)}f(t)^{2}dV_{g(t)} \ \ = \ \ 2\int_{M}f\nabla_{i}f\cdot\nabla^{i}|\nabla u|^{2}dV\\
&=&4\int_{M}f(t)\left\langle du(t)\otimes df(t),\nabla^{2}_{g(t)}u(t)\right\rangle_{g(t)}dV_{g(t)}.
\end{eqnarray*}
Using the contracted Bianchi identities we may simplify the term $\int_{M}\frac{f^{2}\Delta R}{2} dV$ as follows:
\begin{eqnarray*}
&&\int_{N}\frac{f(t)^{2}}{2}\Delta_{g(t)}R_{g(t)} dV_{g(t)} \ \ = \ \
-\frac{1}{2}\int_{M}\nabla_{i}R\cdot\nabla^{i}(f^{2})dV\\
&=&-\int_{M}\nabla_{i}R\cdot f\nabla^{i}f dV \ \ = \ \
-2\int_{M}\nabla^{k}R_{ki}\cdot f\nabla^{i}f dV\\
&=&2\int_{M}R_{ki}\nabla^{k}(f\nabla^{j}f)dV \ \ = \ \ 2\int_{M}R_{ki}
\left(\nabla^{k}f\cdot\nabla^{j}f+f\nabla^{k}\nabla^{j}f\right)dV\\
&=&2\int_{M}\left\langle{\rm Ric}_{g(t)},df(t)\otimes df(t)\right\rangle_{g(t)}dV_{g(t)}\\
&&+2\int_{M}f(t)\left\langle{\rm Ric}_{g(t)},\nabla^{2}_{g(t)}f(t)\right\rangle_{g(t)}dV_{g(t)}.
\end{eqnarray*}
Hence (\ref{8.8}) becomes
\begin{eqnarray*}
\frac{d}{d t}\lambda_{g(t),u(t)}(f(t))&=&\int_{M}\left[2\left\langle{\rm Ric}_{g(t)},df(t)\otimes df(t)
\right\rangle_{g(t)}+\left|{\rm Ric}_{g(t)}\right|^{2}_{g(t)}
f(t)^{2}\right]dV_{g(t)}\\
&&+
 \int_{M}\left[2\left|\Delta_{g(t)}u(t)\right|^{2}_{g(t)}
+4\left|\nabla_{g(t)} u(t)\right|^{4}_{g(t)}\right]f(t)^{2}dV_{g(t)}\\
&&- \ \int_{M}4f(t)^{2}\left\langle du(t)\otimes du(t),{\rm Ric}_{g(t)}\right\rangle_{g(t)}dV_{g(t)}\\
&&- \ \int_{M}4\left|\left\langle \nabla_{g(t)}u(t),\nabla_{g(t)} f(t)\right\rangle_{g(t)}\right|^{2}dV_{g(t)}\\
&=&\int_{M}2\left\langle \mathcal{S}_{g(t)},df(t)\otimes df(t)\right\rangle_{g(t)}dV_{g(t)}\\
&&+\int_{M}f(t)^{2}\left[\left|{\rm Ric}_{g(t)}-2du(t)\otimes du(t)\right|^{2}_{g(t)}+2\left|\Delta_{g(t)}u(t)\right|^{2}_{g(t)}
\right]dV_{g(t)}
\end{eqnarray*}
where by definition $S_{ij}=R_{ij}-2\partial_{i}u\partial_{j}u$.
\end{proof}

In \cite{L}, List proved that the nonnegativity of operator $\mathcal{S}_{g(t)}$ is preserved by the harmonic-Ricci flow, hence

\begin{corollary} \label{c8.3}If ${\rm Ric}_{g(0)}-2du(0)\otimes du(0)\geq0$, then the
eigenvalues of the operator $\Delta_{g(t),u(t)}$ are
nondecreasing under the harmonic-Ricci flow.
\end{corollary}

\begin{remark} \label{r8.4}If we choose $u(t)\equiv0$, then we obtain X. Cao's
result \cite{C1}.
\end{remark}

\section{Another formula for $\frac{d}{dt}\lambda(f(t))$}
\label{section9}

In this section we give another formula for $\frac{d}{dt}\lambda(f(t))$ using the similar method in \cite{Li}. Recall the formula
\begin{eqnarray*}
\frac{d}{d t}\lambda_{g(t),u(t)}(f(t))&=&
\int_{M}2\left\langle\mathcal{S}_{g(t),u(t)},df(t)\otimes df(t)\right\rangle_{g(t)}dV_{g(t)}\\
&&+\int_{M}f(t)^{2}\left[\left|\mathcal{S}_{g(t),u(t)}\right|^{2}_{g(t)}
+2\left|\Delta_{g(t)}u(t)\right|^{2}_{g(t)}\right]dV_{g(t)}.
\end{eqnarray*}
Consider the function $\varphi$ determined by $f^{2}(t)=e^{-\varphi(t)}$. Then we have
\begin{equation*}
df=\frac{-e^{\varphi}d\varphi}{2f}, \ \ \ \frac{\nabla f}{f}=-\frac{\nabla \varphi}{2}, \ \ \ \frac{\Delta f}{f}=-\frac{1}{2}\Delta \varphi
+\frac{1}{4}|\nabla \varphi|^{2}.
\end{equation*}
Hence
\begin{eqnarray*}
2\frac{d}{d t}\lambda_{g(t),u(t)}(f(t))&=&
\int_{M}\left\langle\mathcal{S}_{g(t),u(t)},d\varphi(t)\otimes d\varphi(t)\right\rangle_{g(t)}e^{-\varphi(t)}dV_{g(t)}\\
&&+2\int_{M}\left[\left|\mathcal{S}_{g(t),u(t)}\right|^{2}_{g(t)}
+2\left|\Delta_{g(t)}u(t)\right|^{2}_{g(t)}\right]e^{-\varphi}dV_{g(t)}.
\end{eqnarray*}
Using the integration by parts and contracted Bianchi identities yields
\begin{eqnarray*}
&&\int_{M}\left\langle\mathcal{S}_{g(t),u(t)},d\varphi(t)\otimes d\varphi(t)\right\rangle_{g(t)}e^{-\varphi(t)}dV_{g(t)}\\
&=&\int_{M}S_{ij}\nabla^{i}\varphi\nabla^{j}\varphi e^{-\varphi}dV \ \ = \ \ -\int_{M}S_{ij}\nabla^{j}\varphi\nabla^{i}\left(e^{-\varphi}\right)dV\\
&=&\int_{M}e^{-\varphi}\nabla^{i}\left(S_{ij}\nabla^{j}\varphi\right)dV\\
&=&\int_{M}\nabla^{i}S_{ij}\cdot\nabla^{j}\varphi\cdot e^{-\varphi}dV+\int_{M}S_{ij}\nabla^{i}\nabla^{j}\varphi\cdot e^{-\varphi}dV\\
&=&\int_{M}\nabla^{i}R_{ij}\cdot\nabla^{j}\varphi\cdot e^{-\varphi}dV_{g}+\int_{M}S_{ij}\nabla^{i}\nabla^{j}\varphi\cdot e^{-\varphi}dV\\
&&+\int_{M}\nabla^{i}\left(-2\nabla_{i}u\nabla_{j}u\right)\nabla^{j}\varphi\cdot e^{-\varphi}dV_{g}\\
&=&\frac{1}{2}\int_{M}R\Delta\left(e^{-\varphi}\right)dV
+\int_{M}S_{ij}\nabla^{i}\nabla^{j}\varphi\cdot e^{-\varphi}dV\\
&&-2\int_{M}\left(\nabla^{i}u\nabla_{j}u\right)
\nabla^{i}\nabla^{j}\left(e^{-\varphi}\right)dV.
\end{eqnarray*}
Thus
\begin{eqnarray*}
\int_{M}S_{ij}\nabla^{i}\nabla^{j}\varphi\cdot e^{-\varphi}dV&=&\int_{M}S_{ij}\nabla^{i}\varphi\nabla^{j}\varphi e^{-\varphi}dV-\frac{1}{2}\int_{M}R\Delta\left(e^{-\varphi}\right)dV\\
&&+2\int_{M}\left(\nabla^{i}u\nabla_{j}u\right)\nabla^{i}\nabla^{j}
\left(e^{-\varphi}\right).
\end{eqnarray*}
On the other hand, one gets
\begin{eqnarray*}
&&\int_{M}\left|\nabla^{2}_{g(t)}\varphi(t)\right|^{2}_{g(t)}
e^{-\varphi(t)}dV_{g(t)} \ \ = \ \ \int_{M}\nabla_{i}\nabla_{j}\varphi\nabla^{i}\nabla_{j}\varphi\cdot e^{-\varphi}dV\\
&=&-\int_{M}\nabla_{j}\varphi\cdot\nabla_{i}\nabla^{i}\nabla^{j}
\varphi\cdot e^{-\varphi}dV-\int_{M}\nabla_{j}\varphi\cdot\nabla^{i}\nabla^{j}\varphi\cdot\nabla_{i}\left(e^{-\varphi}\right)dV\\
&=&-\int_{M}\nabla_{j}\varphi\cdot\nabla_{i}\nabla^{j}\nabla^{i}\varphi\cdot e^{-\varphi}dV-\int_{M}\nabla_{j}\varphi\cdot\nabla^{i}\nabla^{j}\varphi\cdot\nabla_{i}\left(e^{-\varphi}\right)dV.
\end{eqnarray*}
Since
\begin{eqnarray*}
\int_{M}\nabla_{j}\varphi\cdot\nabla^{i}\nabla^{j}\varphi\cdot\nabla_{i}
\left(e^{-\varphi}\right)dV&=&-\int_{M}\nabla^{i}\left(\nabla_{j}\varphi\cdot\nabla_{i}
\left(e^{-\varphi}\right)\right)\nabla^{j}\varphi dV\\
&=&-\int_{M}\nabla^{j}\varphi\cdot\nabla^{i}\nabla_{j}\varphi\cdot\nabla_{i}\left(e^{-\varphi}\right)dV
\\
&&- \ \int_{M}\left|\nabla \varphi\right|^{2}\Delta\left(e^{-\varphi}\right)dV
\end{eqnarray*}
which implies
\begin{equation*}
\int_{M}\nabla_{j}\varphi\cdot\nabla^{i}\nabla^{j}\varphi\cdot\nabla_{i}\left(e^{-\varphi}\right)dV
=-\frac{1}{2}\int_{M}\left|\nabla \varphi\right|^{2}
\Delta\left(e^{-\varphi}\right)dV,
\end{equation*}
it follows that
\begin{equation*}
\int_{M}\left|\nabla^{2}\varphi\right|^{2}e^{-\varphi}dV\\
=-\int_{M}\nabla_{j}\varphi\cdot\nabla_{i}\nabla^{j}\nabla^{i}\varphi\cdot e^{-\varphi}dV+\frac{1}{2}\int_{M}\left|\nabla\varphi\right|^{2}\Delta\left(e^{-\varphi}\right)dV.
\end{equation*}
By Ricci identity the term $\nabla^{i}\nabla^{j}\nabla^{i}\varphi$ equals
\begin{eqnarray*}
\nabla_{i}\nabla^{j}\nabla^{i}\varphi&=&g^{jk}g^{il}
\nabla_{i}\nabla_{k}\nabla_{l}\varphi \ \ = \ \ g^{jk}g^{il}\left(\nabla_{k}\nabla_{i}\nabla_{l}\varphi
-R^{p}_{ikl}\nabla_{p}\varphi\right)\\
&=&\nabla^{j}\nabla_{i}\nabla^{i}\varphi-g^{jk}g^{il}R_{iklp}
\nabla^{p}\varphi\\
&=&\nabla^{j}\Delta\varphi+g^{jk}g^{il}R_{ikpl}\nabla^{p}\varphi \ \ = \ \ \nabla^{j}\Delta\varphi+g^{jk}R_{kp}\nabla^{p}\varphi.
\end{eqnarray*}
Hence
\begin{eqnarray*}
&&-\int_{M}\nabla_{j}\varphi\cdot\nabla_{i}\nabla^{j}\nabla^{i}\varphi\cdot e^{-\varphi}dV\\
&=&-\int_{M}\nabla_{i}\varphi\cdot\nabla^{j}\Delta\varphi\cdot
e^{-\varphi}dV-\int_{M}R_{kp}\nabla^{k}\varphi\cdot\nabla^{p}\varphi
e^{-\varphi}dV\\
&=&\int_{M}\nabla^{j}\Delta\varphi\cdot\nabla_{j}\left(e^{-\varphi}\right)
+\int_{M}R_{kp}\nabla^{k}\varphi\cdot\nabla^{p}\left(e^{-\varphi}\right)dV\\
&=&- \ \int_{M}\Delta\varphi\cdot\Delta\left(e^{-\varphi}\right)
-\int_{M}e^{-\varphi}\left(\nabla^{p}R_{kp}\cdot\nabla^{k}\varphi
+R_{kp}\nabla^{p}\nabla^{k}\varphi\right)\\
&=&-\int_{M}\Delta\left(e^{-\varphi}\right)\cdot\Delta\varphi dV+
\frac{1}{2}\int_{M}\nabla_{k}R\cdot\nabla^{k}\left(e^{-\varphi}\right)dV
\\
&&- \ \int_{M}e^{-\varphi}R_{kp}\nabla^{k}\nabla^{p}\varphi dV\\
&=&-\int_{M}\Delta\left(e^{-\varphi}\right)\left(\Delta\varphi
+\frac{1}{2}R\right)-\int_{M}R_{kp}\nabla^{k}\nabla^{p}\varphi\cdot e^{-\varphi}dV.
\end{eqnarray*}
Putting those formulas together, we obtain
\begin{eqnarray*}
&&\int_{M}2S_{ij}\nabla^{i}\nabla^{j}\varphi\cdot e^{-\varphi}dV
+\int_{M}\left|\nabla^{2}\varphi\right|^{2}e^{-\varphi}dV\\
&=&\int_{M}S_{ij}\nabla^{i}\nabla_{j}\varphi\cdot e^{-\varphi}dV
+\int_{M}\left(-2\nabla_{i}u\nabla_{j}u\right)\nabla^{i}\nabla^{j}\varphi\cdot e^{-\varphi}dV\\
&&- \ \int_{M}\Delta\left(e^{-\varphi}\right)\left(\Delta\varphi+\frac{R}{2}
-\frac{1}{2}\left|\nabla\varphi\right|^{2}\right)dV\\
&=&\int_{M}S_{ij}\nabla^{i}\varphi\nabla^{j}
\varphi\cdot e^{-\varphi}dV-\int_{M}\Delta\left(e^{-\varphi}\right)
\left(\Delta\varphi+R-\frac{1}{2}\left|\nabla\varphi
\right|^{2}\right)dV\\
&&+ \ 2\int_{M}\left(\nabla_{i}u\nabla_{j}u\cdot\nabla^{i}\nabla^{j}\left(e^{-\varphi}\right)
-\nabla_{i}u\nabla_{j}u\cdot\nabla^{i}\nabla^{j}\varphi\cdot e^{-\varphi}\right)dV.
\end{eqnarray*}
Since $f$ is an eigenfunction of $\lambda$, it induces $\lambda=-\frac{\Delta f}{f}+\frac{R}{2}-|\nabla
u|^{2}=\frac{1}{2}\Delta\varphi-\frac{1}{4}
|\nabla\varphi|^{2}+\frac{R}{2}-|\nabla u|^{2}$ and therefore
\begin{eqnarray*}
&&\int_{M}2S_{ij}\nabla^{i}\nabla^{j}\varphi\cdot e^{-\varphi}dV
+\int_{M}\left|\nabla^{2}\varphi\right|^{2}e^{-\varphi}dV\\
&=&\int_{M}S_{ij}\nabla^{i}\varphi\nabla^{j}\varphi\cdot e^{-\varphi}dV-2\int_{M}\Delta\left(\left|\nabla u\right|^{2}\right)\cdot e^{-\varphi}dV\\
&&+2\int_{M}\nabla_{i}u\nabla_{j}\left(\nabla^{i}\nabla^{j}\left(e^{-\varphi}\right)
-\nabla^{i}\nabla^{j}\varphi\cdot e^{-\varphi}\right)dV.
\end{eqnarray*}
Plugging into the expression of $\frac{d}{dt}\lambda(f(t))$ yields
\begin{eqnarray*}
&&2\frac{d}{d t}\lambda_{g(t),u(t)}(f(t)) \ \ = \ \ \int_{M}S_{ij}\nabla^{i}\varphi\nabla^{j}\varphi\cdot e^{-\varphi}dV
+\int_{M}\left|\mathcal{S}\right|^{2}e^{-\varphi}dV\\
&&+ \ \int_{M}\left|\mathcal{S}\right|^{2}e^{-\varphi}dV
+4\int_{M}\left|\Delta u\right|^{2}e^{-\varphi}dV\\
&=&\int_{M}\left|\mathcal{S}_{g(t),u(t)}+\nabla^{2}_{g(t)}\varphi(t)
\right|^{2}_{g(t)}e^{-\varphi(t)}dV_{g(t)}+\int_{M}
\left|\mathcal{S}_{g(t),u(t)}\right|^{2}_{g(t)}e^{-\varphi(t)}dV_{g(t)}\\
&&+ \ 4\int_{M}\left|\Delta_{g(t)}u(t)\right|^{2}_{g(t)}e^{-\varphi(t)}
dV_{g(t)}+2\int_{M}\Delta_{g(t)}\left|\nabla_{g(t)} u(t)\right|^{2}_{g(t)}e^{-\varphi(t)}dV_{g(t)}\\
&&+ \ 2\int_{M}\nabla_{i}u\nabla_{j}u\left[-\nabla^{i}\nabla^{j}
\left(e^{-\varphi}\right)+\nabla^{i}\nabla^{j}\varphi\cdot e^{-\varphi}\right]dV
\end{eqnarray*}
On the other hand,
\begin{eqnarray*}
I&:=&\int_{M}\left(\nabla_{i}u\nabla_{j}u\cdot\nabla^{i}\nabla^{j}\varphi\right)e^{-\varphi}dV\\
&=&-\int_{M}\nabla^{i}\left(\nabla_{i}u\nabla_{j}u\cdot e^{-\varphi}\right)\nabla^{j}\varphi dV\\
&=&-\int_{M}\nabla^{j}\varphi\left(\Delta u\cdot\nabla_{j}
u\cdot e^{-\varphi}+\nabla_{i}u\nabla^{i}\nabla_{j}u\cdot e^{-\varphi}
-\nabla_{i}u\nabla_{j}u\nabla^{i}\varphi\cdot e^{-\varphi}\right)dV\\
&=&-\int_{M}\nabla_{j}u\nabla^{j}\varphi\Delta u\cdot e^{-\varphi}dV
-\int_{M}\nabla_{i}u\nabla^{j}\varphi\nabla^{i}\nabla_{j}u
\cdot e^{-\varphi}dV\\
&&+\int_{M}\left|\langle du,d\varphi\rangle\right|^{2}e^{-\varphi}dV
\end{eqnarray*}
and
\begin{eqnarray*}
II&:=&\int_{M}\nabla_{i}u\nabla_{j}u\nabla^{i}\nabla^{j}
\left(e^{-\varphi}\right)dV \ \ = \ \ \int_{M}\nabla^{i}\nabla^{j}\left(\nabla_{i}u\nabla_{j}u\right)e^{-\varphi}dV\\
&=&\int_{M}\nabla^{i}\left(\nabla^{j}\nabla_{i}u\cdot\nabla_{j}u
+\nabla_{i}u\Delta u\right)e^{-\varphi}dV\\
&=&\int_{M}\left(\Delta\nabla^{i}u\cdot\nabla_{i}u+
\nabla^{i}\Delta u\cdot\nabla_{i}u+\left|\nabla^{2}u\right|^{2}
+\left|\Delta u\right|^{2}\right)e^{-\varphi}dV
\end{eqnarray*}
and
\begin{eqnarray*}
III&:=&\int_{M}\Delta\left(\left|\nabla u\right|^{2}\right)
e^{-\varphi}dV \ \ = \ \ 2\int_{M}\nabla^{i}\left(\nabla_{i}\nabla_{j}u\cdot\nabla^{j}
u\right)e^{-\varphi}dV\\
&=&2\int_{M}\left(\Delta\nabla_{j}u\cdot\nabla^{j}u+\left|\nabla^{2}u\right|^{2}\right)e^{-\varphi}dV.
\end{eqnarray*}
If we set
\begin{equation*}
B:=2\left(III+I-II\right)
\end{equation*}
then
\begin{eqnarray*}
\frac{B}{2}&=&\int_{M}\left[\Delta\nabla_{i}u\cdot\nabla^{i}u
-\nabla_{i}\Delta u\cdot\nabla^{i}u+\left|\nabla^{2}u\right|^{2}
-\left|\Delta u\right|^{2}+
\left|\langle du,d\varphi\rangle\right|^{2}\right.\\
&&\left.-\nabla_{i}u\cdot\nabla^{i}\varphi\cdot\Delta u
-\nabla_{i}u\cdot\nabla^{j}\varphi\cdot\nabla^{i}\nabla_{j}u\right]e^{-\varphi}dV\\
&=&\int_{M}\left(R_{ij}\nabla^{i}u\nabla^{j}u+\left|\nabla^{2}u\right|^{2}
-\left|\Delta u\right|^{2}+\left|\langle du,d\varphi\rangle\right|^{2}\right.\\
&&\left.-\nabla_{i}u\cdot\nabla^{i}\varphi\cdot\Delta u-\nabla_{i}u
\cdot\nabla^{j}\varphi\cdot\nabla^{i}\nabla_{j}u\right)
e^{-\varphi}dV.
\end{eqnarray*}
On the other hand,
\begin{eqnarray*}
-\int_{M}\nabla_{i}u\cdot\nabla^{i}\varphi\cdot\Delta u\cdot e^{-\varphi}dV
&=&\int_{M}\left(\nabla_{i}u\cdot\Delta u\right)\nabla^{i}\left(e^{-\varphi}\right)dV\\
&=&-\int_{M}\nabla^{i}\left(\nabla_{i}u\cdot\Delta u\right)e^{-\varphi}dV\\
&=&\int_{M}\left(-\left|\Delta u\right|^{2}-\nabla_{i}u\cdot\nabla^{i}\Delta u\right)e^{-\varphi}dV
\end{eqnarray*}
and
\begin{eqnarray*}
-\int_{M}\nabla_{i}u\nabla^{j}\varphi\nabla^{i}\nabla_{j}u\cdot
e^{-\varphi}dV
&=&\int_{M}\nabla_{i}u\nabla^{i}\nabla_{j}u\nabla^{j}
\left(e^{-\varphi}\right)dV\\
&=&-\int_{M}\nabla^{j}\left(\nabla_{i}u\nabla^{i}\nabla_{j}u\right)e^{-\varphi}dV\\
&=&\int_{M}\left(-\left|\nabla^{2}u\right|^{2}
-\nabla_{i}u\Delta\nabla^{i}u\right)e^{-\varphi}dV.
\end{eqnarray*}
Therefore
\begin{equation}
\frac{B}{2}=\int_{M}\left[-2\left|\Delta u\right|^{2}
+\left|\langle du,d\varphi\rangle\right|^{2}-2\langle\nabla u,
\nabla\Delta u\rangle\right]e^{-\varphi}dV.\label{9.1}
\end{equation}
By definition,
\begin{equation*}
\Delta\left(\left|\nabla u\right|^{2}\right)=\Delta\left(\nabla^{i}u
\cdot\nabla_{i}u\right)=2\nabla^{i}u\cdot\Delta\nabla_{i}u+2\left|\nabla^{2}u\right|^{2}.
\end{equation*}
So
\begin{eqnarray*}
\Delta\left|\nabla u\right|^{2}&=&2\left|\nabla^{2}u\right|^{2}
+2\left(\nabla_{i}\Delta u+R_{ij}\nabla^{j}u\right)\nabla^{i}u\\
&=&2\left|\nabla^{2}u\right|^{2}+2R_{ij}\nabla^{i}u\cdot\nabla^{j}u+
2\langle\nabla u,\nabla\Delta u\rangle.
\end{eqnarray*}.
Pugging it into (\ref{9.1}) yields
\begin{equation*}
\frac{B}{2}=\int_{M}\left[-2\left|\Delta u\right|^{2}
+\left|\langle du,d\varphi\rangle\right|^{2}+
2\left|\nabla^{2}u\right|^{2}-\Delta\left|\nabla u\right|^{2}+2R_{ij}\nabla^{i}u\nabla^{j}u\right]e^{-\varphi}dV.
\end{equation*}
Since
\begin{eqnarray*}
2R_{ij}\nabla^{i}u\nabla^{j}u&=&2\left(S_{ij}+2\nabla_{i}u\nabla_{j}u\right)
\nabla^{i}u\nabla^{j}u\\
&=&2S_{ij}\nabla^{i}u\nabla^{j}u+4\left|\nabla u\right|^{4}\\
&=&\frac{1}{4}\left|\mathcal{S}+4du\otimes du\right|^{2}-\frac{1}{4}\left|\mathcal{S}\right|^{2},
\end{eqnarray*}
it follows that
\begin{eqnarray*}
\frac{B}{2}&=&III+I-II\\
&=&\int_{M}\left[\left|\langle du,d\varphi\rangle\right|^{2}
-2|\Delta u|^{2}-\frac{1}{4}\left|\mathcal{S}\right|^{2}+2\left|\nabla^{2}u\right|^{2}+\frac{1}{4}\left|\mathcal{S}+4du\otimes du\right|^{2}\right]e^{-\varphi}dV\\
&&- \ III.
\end{eqnarray*}
Hence
\begin{eqnarray*}
B&=&\int_{M}\left[-4\left|\Delta u\right|^{2}+2\left|\langle du,
d\varphi\rangle\right|^{2}-\frac{1}{2}\left|\mathcal{S}\right|^{2}\right.\\
&& \ +\left.4\left|\nabla^{2}u\right|^{2}+\frac{1}{2}\left|\mathcal{S}
+4du\otimes du\right|^{2}\right]e^{-\varphi}dV-2III.
\end{eqnarray*}

\begin{theorem} \label{t9.1}Suppose that $(g(t),u(t))$ is a solution of the harmonic-Ricci flow on a compact Riemannian manifold $M$ and $f(t)$ is an eigenfunction of $\Delta_{g(t),u(t)}$, i.e., $\Delta_{g(t),u(t)}f(t)=\lambda(t)f(t)$(where
$\lambda(t)$ is only a function of time $t$), with the normalized condition $\int_{M}f(t)^{2}dV_{g(t)}=1$. Then we have
\begin{eqnarray}
&&\frac{d}{dt}\lambda(t) \ \ = \ \ \frac{d}{d t}\lambda_{g(t),u(t)}(f(t))
\label{9.2}\\
&=&\frac{1}{2}\int_{M}\left|\mathcal{S}_{g(t),u(t)}
+\nabla^{2}_{g(t)}\varphi(t)\right|^{2}_{g(t)}
e^{-\varphi(t)}dV_{g(t)}\nonumber\\
&&+ \ \frac{1}{4}\int_{M}\left|\mathcal{S}_{g(t),u(t)}\right|^{2}_{g(t)}
e^{-\varphi(t)}dV_{g(t)}\nonumber\\
&&+ \ \int_{M}\left|\langle du(t),d\varphi(t)\rangle_{g(t)}\right|^{2}e^{-\varphi(t)}dV_{g(t)}
+2\int_{M}\left|\nabla^{2}_{g(t)}u(t)\right|^{2}_{g(t)}
e^{-\varphi(t)}dV_{g(t)}\nonumber\\
&&+ \ \frac{1}{4}\int_{M}\left|\mathcal{S}_{g(t),u(t)}+4du(t)\otimes du(t)\right|^{2}_{g(t)}e^{-\varphi(t)}dV_{g(t)}\nonumber\\
&&- \ \int_{M}\Delta_{g(t)}\left(\left|\nabla_{g(t)} u(t)\right|^{2}_{g(t)}\right)e^{-\varphi(t)}dV_{g(t)}.\nonumber
\end{eqnarray}

\end{theorem}

\begin{remark} \label{r5.2} When $u\equiv0$, (\ref{9.2}) reduces to
J. Li's formula \cite{Li}.
\end{remark}

\section{The first variation of expander and shrinker entropies}
\label{section10}

Suppose that $M$ is a closed manifold of dimension $n$. We define
\begin{equation*}
\mathcal{W}_{\pm}: \odot^{2}_{+}(M)\times C^{\infty}(M)\times C^{\infty}(M)\times\mathbb{R}^{+}\longrightarrow\mathbb{R}, \ \ \ (g,u,f,\tau)
\longmapsto\mathcal{W}_{\pm}(g,u,f,\tau)
\end{equation*}
where
\begin{equation}
\mathcal{W}_{\pm}(g,u,f,\tau):=\int_{M}\left[\tau\left(S_{g,u}
+\left|\nabla_{g} f\right|^{2}_{g}\right)\mp f\pm n\right]\frac{e^{-f}}{(4\pi\tau)^{n/2}}dV_{g}.\label{10.1}
\end{equation}
Set
\begin{eqnarray*}
\mu_{\pm}(g,u,\tau)&:=&\inf\left\{\mathcal{W}_{\pm}(g,u,f,\tau)\Big| f\in C^{\infty}(M), \ \ \ \int_{M}\frac{e^{-f}}{(4\pi\tau)^{n/2}}dV_{g}=1\right\},\\
\nu_{\pm}(g,u)&:=&\sup\{\mu_{\pm}(g,u,\tau)|\tau>0\}.
\end{eqnarray*}

\begin{lemma} \label{l10.1}Suppose $\nu_{\pm}(g,u)=\mathcal{W}_{\pm}(g,u,f_{\pm},\tau_{\pm})$ for some functions $f_{\pm}$ and constants $\tau_{\pm}$ satisfying
\begin{equation*}
\int_{M}\frac{e^{-f_{\pm}}}{(4\pi\tau_{\pm})^{n/2}}dV_{g}=1, \ \ \ \tau_{\pm}>0,
\end{equation*}
then we must have
\begin{eqnarray*}
\tau_{\pm}\left(-2\Delta_{g}f_{\pm}+\left|\nabla_{g} f_{\pm}\right|^{2}_{g}
-S_{g,u}\right)\pm f_{\pm}\mp n+\nu_{\pm}(g,u)&=&0,\\
\int_{M}\frac{f_{\pm}e^{-f_{\pm}}}{(4\pi\tau)^{n/2}}dV_{g}&=&\frac{n}{2}\mp\nu_{\pm}(g,u).
\end{eqnarray*}
\end{lemma}

\begin{proof} Since $g$ and $u$ are fixed, we consider the corresponding
Lagrangian multiplier function
\begin{equation*}
\mathfrak{L}_{\pm}(f,\tau;\lambda):=\mathcal{W}_{\pm}(g,u,f,\tau)-\lambda\left(\int_{M}\frac{e^{-f}}{(4\pi\tau)^{n/2}}dV_{g}
-1\right).
\end{equation*}
Then the variation of $\mathfrak{L}_{\pm}$ in $f$ direction is
\begin{eqnarray*}
\delta_{f}\mathfrak{L}_{\pm}(f,\tau;\lambda)&=&\int_{M}\left[2\tau\nabla^{i}f\nabla_{i}(\delta f)\mp\delta f+\lambda\delta f\right]\frac{e^{-f}}{(4\pi\tau)^{n/2}}dV_{g}\\
&&- \ \int_{M}\left[\tau\left(S_{g,u}+\left|\nabla_{g} f\right|^{2}_{g}\right)
\mp f\pm n\right]\delta f\frac{e^{-f}}{(4\pi\tau)^{n/2}}dV_{g}.
\end{eqnarray*}
By the divergence theorem, we calculate
\begin{eqnarray*}
\int_{M}\nabla^{i}f\cdot\nabla_{i}(\delta f)\frac{e^{-f}}{(4\pi\tau)^{n/2}}dV_{g}&=&-\int_{M}\nabla_{i}\left(\nabla^{i}f\frac{e^{-f}}{(4\pi\tau)^{n/2}}\right)\delta f dV_{g}\\
&=&-\int_{M}\left(\Delta_{g}f-\left|\nabla_{g} f\right|^{2}_{g}\right)\delta f\frac{e^{-f}}{(4\pi\tau)^{n/2}}dV_{g}.
\end{eqnarray*}
Hence
\begin{equation*}
\delta_{f}\mathfrak{L}_{\pm}(f,\tau;\lambda)=\int_{M}
\left[\tau\left(-2\Delta_{g}f+\left|\nabla_{g} f\right|^{2}_{g}-S_{g,u}\right)\pm f\mp n\mp 1+\lambda\right]\delta f\frac{e^{-f}}{(4\pi\tau)^{n/2}}dV.
\end{equation*}
This implies that
\begin{equation*}
\tau_{\pm}\left(-2\Delta_{g}f_{\pm}+\left|\nabla_{g} f_{\pm}\right|^{2}_{g}-S_{g,u}\right)\pm f_{\pm}\mp n\mp 1+\lambda_{\pm}=0.
\end{equation*}
Since $f_{\pm}$ satisfies the normalized condition, it follows that
\begin{equation*}
0=\lambda_{\pm}\mp1+\int_{M}\left[\tau_{\pm}\left(-2\Delta_{g}
f_{\pm}+\left|\nabla_{g} f_{\pm}\right|^{2}_{g}-S_{g,u}\right)\pm f_{\pm}\mp n\right]\frac{e^{-f_{\pm}}}{(4\pi\tau_{\pm})^{n/2}}dV_{g}.
\end{equation*}
From the identity
\begin{equation*}
\int_{M}\Delta_{g}f\frac{e^{-f}}{(4\pi\tau)^{n/2}}dV_{g}=
\int_{M}\left|\nabla_{g} f\right|^{2}_{g}\frac{e^{-f}}{(4\pi\tau)^{n/2}}dV_{g}
\end{equation*}
and the definition (\ref{10.1}), we obtain
\begin{equation*}
\nu_{\pm}(g,u)=\mathcal{W}_{\pm}(g,u,f_{\pm},\tau_{\pm})=\lambda_{\pm}\mp1,
\end{equation*}
and consequently,
\begin{equation*}
\tau_{\pm}\left(-2\Delta_{g}f_{\pm}+\left|\nabla_{g} f_{\pm}\right|^{2}_{g}
-S_{g,u}\right)\pm f_{\pm}\mp n+\nu_{\pm}(g,u)=0.
\end{equation*}
The variation of $\mathfrak{L}_{\pm}$ with respect to $\tau$ indicates
\begin{eqnarray*}
\delta_{\tau}\mathfrak{L}_{\pm}(f,\tau;\lambda)&=&\int_{M}\delta\tau\left(S_{g,u}
+\left|\nabla_{g} f\right|^{2}_{g}\right)\frac{e^{-f}}{(4\pi\tau)^{n/2}}dV_{g}
\\
&&- \ \lambda\int_{M}\left(-\frac{n}{2}\frac{\delta\tau}{\tau}\right)\frac{e^{-f}}{(4\pi\tau)^{n/2}}dV_{g}.
\\
&&+ \ \int_{M}\left(-\frac{n}{2}\frac{\delta\tau}{\tau}\right)\left[\tau\left(S_{g,u}
+\left|\nabla_{g} f\right|^{2}_{g}\right)\mp f\pm n\right]\frac{e^{-f}}{(4\pi\tau)^{n/2}}dV_{g}\\
&=&\int_{M}\delta\tau\left[\left(1-\frac{n}{2}\right)
\left(S_{g,u}+\left|\nabla_{g} f\right|^{2}_{g}\right)+\frac{n}{2\tau}\left(\lambda\pm f\mp n\right)\right]\frac{e^{-f}dV_{g}}{(4\pi\tau)^{n/2}}.
\end{eqnarray*}
Using the first proved equation we have
\begin{eqnarray*}
0&=&\int_{M}\left[\left(\nu_{\pm}(g,u)\pm f_{\pm}\mp n\right)\left(1-\frac{n}{2}\right)+\frac{n}{2}\left(\nu_{\pm}(g,u)\pm f_{\pm}\mp n\pm1\right)\right]\frac{e^{-f_{\pm}}dV_{g}}{(4\pi\tau_{\pm})^{n/2}}
\\
&=&\int_{M}\left(\nu_{\pm}\pm f_{\pm}\mp\frac{n}{2}\right)\frac{e^{-f_{\pm}}}{(4\pi\tau_{\pm})^{n/2}}
dV_{g}
\end{eqnarray*}
and therefore we obtain the second one.
\end{proof}

For a symmetric $2$-tensor $h=(h_{ij})\in\odot^{2}(M)$, we set
\begin{equation*}
g(s):=g+sh
\end{equation*}
Then the variation of $g(s)$ is
\begin{equation}
\frac{\partial}{\partial s}\Big|_{s=0}R_{g(s)}
=-h^{ij}R_{ij}+\nabla^{i}\nabla^{j}h_{ij}-\Delta_{g}\left({\rm tr}_{g}h\right).\label{10.2}
\end{equation}

\begin{theorem} \label{t10.2}Suppose that $(M,g)$ is a compact Riemannian manifold and $u$ a smooth function on $M$. Let $h$ be any
symmetric covariant $2$-tensor on $M$ and set $g(s):=g+sh$. Let $v$ be any smooth function on $M$ and $u(s):=u+sv$. If $\nu_{\pm}(g(s),u(s))=\mathcal{W}_{\pm}(g(s),u(s),f_{\pm}(s),\tau_{\pm}(s))$ for some smooth functions $f_{\pm}(s)$ with $\int_{M}e^{-f_{\pm}(s)}dV/(4\pi\tau_{\pm}(s))^{n/2}=1$ and constants $\tau_{\pm}(s)>0$, then
\begin{eqnarray*}
&&\frac{d}{ds}\Big|_{s=0}\nu_{\pm}(g(s),u(s))\\
&=&-\tau_{\pm}\int_{M}
\left(\left\langle h,\mathcal{S}_{g,u}\right\rangle_{g}+\left\langle h,\nabla^{2}_{g}f\right\rangle_{g}\pm\frac{1}{2\tau_{\pm}}{\rm tr}_{g}h\right)\frac{e^{-f_{\pm}}}{(4\pi\tau_{\pm})^{n/2}}dV_{g}\\
&&+4\tau_{\pm}\int_{M}v\left(\Delta_{g} u-\langle du,df_{\pm}\rangle_{g}\right)\frac{e^{-f_{\pm}}}{(4\pi\tau_{\pm})^{n/2}}dV_{g},
\end{eqnarray*}
where $f_{\pm}:=f_{\pm}(0)$ and $\tau_{\pm}:=\tau_{\pm}(0)$. In particular, the critical points of $\nu_{\pm}(\cdot,\cdot)$ satisfy
\begin{equation*}
\mathcal{S}_{g,u}+\nabla^{2}_{g}f\pm\frac{1}{2\tau_{\pm}}g=0, \ \ \ \Delta_{g}u=\left\langle du,df_{\pm}\right\rangle_{g}.
\end{equation*}
Consequently, if $\mathcal{W}_{\pm}(g,u,f,\tau)$ and
$\nu_{\pm}(g,u)$ achieve their minimums, then $(M,g)$ is a gradient expanding and shrinker harmonic-Ricci soliton according to the sign.
\end{theorem}

\begin{proof} By definition, one has
\begin{eqnarray*}
&&\frac{d}{d s}\nu_{\pm}(g(s),u(s)) \ \ = \ \ \frac{d}{d s}\mathcal{W}_{\pm}(g(s),u(s),f_{\pm}(s),\tau_{\pm}(s))\\
&=&\int_{M}\left[\frac{\partial}{\partial s}\tau_{\pm}(s)\left(S_{g(s),u(s)}+\left|\nabla_{g(s)} f_{\pm}(s)\right|^{2}_{g(s)}\right)\right]\frac{e^{-f_{\pm}(s)}}{(4\pi\tau_{\pm}(s))^{n/2}}dV_{g(s)}\\
&&+\int_{M}\left[\tau_{\pm}(s)\frac{\partial}{\partial s}\left(S_{g(s),u(s)}+\left|\nabla_{g(s)} f_{\pm}(s)\right|^{2}_{g(s)}\right)\mp\frac{\partial}{\partial s}f_{\pm}(s)\right]\frac{e^{-f_{\pm}(s)}}{(4\pi\tau_{\pm}(s))^{n/2}}dV_{g(s)}\\
&&+\int_{M}\left[\tau_{\pm}(s)\left(S_{g(s),u(s)}
+\left|\nabla_{g(s)} f_{\pm}(s)\right|^{2}_{g(s)}\right)\mp f_{\pm}(s)\pm n\right]\\
&&\cdot\frac{\partial}{\partial s}\left(\frac{e^{-f_{\pm}(s)}}{(4\pi\tau_{\pm}(s))^{n/2}} dV_{g(s)}\right).
\end{eqnarray*}
Since
\begin{eqnarray*}
\frac{\partial}{\partial s}S_{g(s),u(s)}&=&\frac{\partial}{\partial s}R_{g(s)}
-2\frac{\partial}{\partial s}\left|\nabla_{g(s)} u(s)\right|^{2}_{g(s)}\\
&=&\frac{\partial}{\partial s}R_{g(s)}-2\left(\frac{\partial}{\partial s}g^{ij}\right)\nabla_{i}u\nabla_{j}u-4g^{ij}\frac{\partial}{\partial s}\nabla_{i}u\cdot\nabla_{j}u\\
&=&\frac{\partial}{\partial s}R_{g(s)}-2\left(-g^{ip}g^{jq}h_{pq}\right)\nabla_{i}u\nabla_{j}u-4g^{ij}\nabla_{i}\left(\frac{\partial}{\partial s}u\right)\nabla_{j}u\\
&=&\frac{\partial}{\partial s}R_{g(s)}+2h_{pq}\nabla^{p}u\nabla^{q}u-4\nabla_{i}\left(\frac{\partial}{\partial t}u\right)\nabla^{i}u
\end{eqnarray*}
and
\begin{eqnarray*}
&&\frac{\partial}{\partial s}\left(\frac{e^{-f_{\pm}(s)}}{(4\pi\tau_{\pm}(s))^{n/2}}dV_{g(s)}\right)\\
&=&\left(-\frac{\partial}{\partial s}f_{\pm}(s)-\frac{n}{2\tau_{\pm}(s)}
\frac{\partial}{\partial s}\tau_{\pm}(s)\right)\frac{e^{-f_{\pm}(s)}}{(4\pi\tau_{\pm}(s))^{n/2}}dV_{g(s)}\\
&&+\frac{e^{-f_{\pm}(s)}}{(4\pi\tau_{\pm}(s))^{n/2}}\frac{\partial}{\partial s}dV_{g(s)}\\
&=&\left(-\frac{\partial}{\partial s}f_{\pm}(s)-\frac{n}{2\tau_{\pm}(s)}\frac{\partial}{\partial s}\tau_{\pm}(s)+\frac{1}{2}{\rm tr}_{g}h\right)\frac{e^{-f_{\pm}(s)}}{(4\pi\tau_{\pm}(s))^{n/2}}dV_{g(s)},
\end{eqnarray*}
it follows that
\begin{eqnarray*}
&&\frac{d}{d s}\nu_{\pm}(g(s),u(s))\\
&=&\int_{M}\frac{\partial}{\partial s}\tau_{\pm}(s)\left(S_{g(s),u(s)}+\left|\nabla_{g(s)} f_{\pm}(s)\right|^{2}_{g(s)}\right)\frac{e^{-f_{\pm}(s)}}{(4\pi\tau_{\pm}(s))^{n/2}}dV_{g(s)}\\
&&+\int_{M}\left[\tau_{\pm}(s)\left(\frac{\partial}{\partial s}R_{g(s)}+2h_{pq}\nabla^{p}u\nabla^{q}u-4\nabla_{i}\left(\frac{\partial}{\partial s}u\right)\nabla^{i}u\right.\right.\\
&&\left.\left.-h_{pq}\nabla^{p}f\nabla^{q}f+2\nabla_{i}\left(\frac{\partial}{\partial s}f\right)\nabla^{i}f\right)\mp\frac{\partial}{\partial s}f_{\pm}(s)\right]\frac{e^{-f_{\pm}(s)}}{(4\pi\tau_{\pm}(s))^{n/2}}dV_{g(s)}\\
&&+\int_{M}\left(-\frac{\partial}{\partial s}f_{\pm}(s)-\frac{n}{2\tau_{\pm}(s)}\frac{\partial}{\partial s}\tau_{\pm}(s)+\frac{1}{2}{\rm tr}_{g}h\right)\cdot\\
&&\left[\tau_{\pm}(s)\left(S_{g(s),u(s)}+\left|\nabla_{g(s)} f_{\pm}(s)\right|^{2}_{g(s)}\right)\mp f_{\pm}(s)\pm n\right]\frac{e^{-f_{\pm}(s)}}{(4\pi\tau_{\pm}(s))^{n/2}}dV_{g(s)}.
\end{eqnarray*}
Since
\begin{eqnarray*}
\int_{M}\Delta_{g}{\rm tr}_{g}h\cdot e^{-f}dV_{g}&=&\int_{M}{\rm tr}_{g}h\cdot\Delta_{g}\left(e^{-f}\right)dV_{g} \\
&=&\int_{M}{\rm tr}_{g}h\left(-\Delta_{g}f+\left|\nabla_{g} f\right|^{2}_{g}\right)e^{-f}dV_{g},\\
\int_{M}\nabla^{i}\nabla^{j}h_{ij}\cdot e^{-f}dV_{g}&=&\int_{M}h_{ij}\nabla^{i}\nabla^{j}\left(e^{-f}\right)dV\\
&=&\int_{M}h_{ij}\left(-\nabla^{i}\nabla^{j}f+\nabla^{i}f\nabla^{j}f\right)e^{-f}
dV_{g},\\
\int_{M}\nabla_{i}\left(\frac{\partial}{\partial s}f\right)\nabla^{i}fe^{-f}dV_{g}
&=&\int_{M}-\frac{\partial}{\partial s}f\left(\Delta_{g}f-\left|\nabla_{g} f\right|^{2}_{g}\right)e^{-f}dV_{g},\\
\int_{M}\Delta_{g}\left(e^{-f}\right)dV_{g}&=&\int_{M}
\left(-\Delta_{g}f+\left|\nabla_{g} f\right|^{2}_{g}\right)e^{-f}dV_{g},
\end{eqnarray*}
and using Lemma \ref{l10.1}, we obtain
\begin{eqnarray*}
&&\frac{d}{ds}\Big|_{s=0}\nu_{\pm}(g(s),u(s))\\
&=&\int_{M}\frac{\partial}{\partial s}\Big|_{s=0}\tau_{\pm}(s)\left(S_{g,u}+\left|\nabla_{g} f\right|^{2}_{g}\right)
\frac{e^{-f_{\pm}}}{(4\pi\tau_{\pm})^{n/2}}dV_{g}\\
&&+ \ \int_{M}\bigg[\tau_{\pm}\bigg(-h^{ij}R_{ij}+\nabla^{i}\nabla_{j}h_{ij}
-\Delta_{g}\left({\rm tr}_{g}h\right)+2h_{pq}\nabla^{p}u\nabla^{q}u\\
&& \ -4\nabla_{i}v\nabla^{i}u-h_{pq}\nabla^{p}f\nabla^{q}f
+2\nabla_{i}\left(\frac{\partial}{\partial s}\Big|_{s=0}f(s)\right)\nabla^{i}f\bigg)\mp\frac{\partial}{\partial s}\Big|_{s=0}f(s)\bigg]\\
&& \ \frac{e^{-f_{\pm}}}{(4\pi\tau_{\pm})^{n/2}}dV_{g}\\
&&+\int_{M}\left(-\frac{\partial}{\partial s}\Big|_{s=0}f_{\pm}(s)
-\frac{n}{2\tau_{\pm}}(s)\frac{\partial}{\partial s}\Big|_{s=0}\tau_{\pm}(s)
+\frac{1}{2}{\rm tr}_{g}h\right)\\
&&\cdot\left[\tau_{\pm}\left(S_{g,u}+\left|\nabla_{g} f_{\pm}\right|^{2}_{g}\right)
\mp f_{\pm}\pm n\right]\frac{e^{-f_{\pm}}}{(4\pi\tau_{\pm})^{n/2}}dV_{g}.
\end{eqnarray*}
If we denote by $B$ the last term while $A$ the rest terms, then
\begin{eqnarray*}
A&=&\int_{M}\left[\frac{\partial}{\partial s}\Big|_{s=0}\tau_{\pm}(s)\left(\left|\nabla_{g} f_{\pm}\right|^{2}_{g}+S_{g,u}\right)\right.\\
&& \ -\left.\tau_{\pm}\left(h^{ij}\nabla_{i}\nabla_{j}f_{\pm}
+h^{ij}S_{ij}+4\nabla_{i}v\cdot\nabla^{i}u\right)\mp\frac{\partial}{\partial s}f_{\pm}\right]\frac{e^{-f_{\pm}}}{(4\pi\tau_{\pm})^{n/2}}dV_{g}\\
&&+ \ \int_{M}\tau_{\pm}\left(\Delta_{g}f_{\pm}
-\left|\nabla_{g} f_{\pm}\right|^{2}_{g}\right)\left({\rm tr}_{g}h-2\frac{\partial}{\partial s}\Big|_{s=0}f(s)\right)\frac{e^{-f_{\pm}}}{(4\pi\tau_{\pm})^{n/2}}dV_{g}.
\end{eqnarray*}
The normalized condition
\begin{equation*}
1=\int_{M}\frac{e^{-f_{\pm}(s)}}{(4\pi\tau_{\pm}(s))^{n/2}}dV_{g}
\end{equation*}
implies
\begin{equation*}
0=\int_{M}\left(-\frac{\partial}{\partial s}\Big|_{s=0}f_{\pm}(s)-\frac{n}{2\tau_{\pm}}\frac{\partial}{\partial s}\Big|_{s=0}\tau_{\pm}(s)+\frac{1}{2}{\rm tr}_{g}h\right)\frac{e^{-f_{\pm}(s)}}{(4\pi\tau_{\pm}(s))^{n/2}}dV_{g}.
\end{equation*}
Lemma \ref{l10.1} concludes that
\begin{equation*}
\tau_{\pm}S_{g,u}-\tau_{\pm}\left(\left|\nabla_{g} f_{\pm}\right|^{2}_{g}
-2\Delta_{g}f_{\pm}\right)=\pm f_{\pm}\mp n+\nu_{\pm}(g,u)
\end{equation*}
therefore
\begin{equation*}
\tau_{\pm}\left(S_{g,u}+\left|\nabla_{g} f_{\pm}\right|^{2}_{g}\right)\mp f_{\pm}\pm n=2\tau_{\pm}\left(\left|\nabla_{g} f_{\pm}\right|^{2}_{g}-\Delta_{g}f_{\pm}\right)+\nu_{\pm}(g,u)
\end{equation*}
Plugging it into the definition of $B$ yields
\begin{eqnarray*}
B&=&\int_{M}\left(-\frac{\partial}{\partial s}\Big|_{s=0}f_{\pm}(s)-\frac{n}{2\tau_{\pm}}\frac{\partial}{\partial s}\Big|_{s=0}\tau_{\pm}(s)+\frac{1}{2}{\rm tr}_{g}h\right)\\
&&\cdot\left[2\tau_{\pm}\left(\left|\nabla_{g} f_{\pm}\right|^{2}_{g}-\Delta_{g}f_{\pm}\right)+\nu_{\pm}(g,u)\right]
\frac{e^{-f_{\pm}}}{(4\pi\tau_{\pm})^{n/2}}dV_{g}\\
&=&\int_{M}\left(-\frac{\partial}{\partial s}\Big|_{s=0}f_{\pm}(s)-\frac{n}{2\tau_{\pm}}\frac{\partial}{\partial s}\Big|_{s=0}\tau_{\pm}(s)+\frac{1}{2}{\rm tr}_{g}h\right)\\
&&\cdot\left[2\tau_{\pm}\left(\left|\nabla_{g} f_{\pm}\right|^{2}_{g}-\Delta_{g}f_{\pm}\right)\right]
\frac{e^{-f_{\pm}}}{(4\pi\tau_{\pm})^{n/2}}dV_{g}\\
&=&\int_{M}\left(-\frac{\partial}{\partial s}\Big|_{s=0}f_{\pm}(s)+\frac{1}{2}{\rm tr}_{g}h\right)2\tau_{\pm}\left(\left|\nabla_{g} f_{\pm}\right|^{2}_{g}-\Delta_{g}f_{\pm}\right)
\frac{e^{-f_{\pm}}}{(4\pi\tau_{\pm})^{n/2}}dV_{g}
\end{eqnarray*}
where we use the fact that $\int_{M}\Delta_{g}(e^{-f})dV_{g}=0$. Hence $B$ cancels with the last term in $A$. Therefore, the above variation equals
\begin{eqnarray*}
&&\frac{d}{ds}\Big|_{s=0}\nu_{\pm}(g(s),u(s))\\
&=&\int_{M}\bigg[\frac{\partial}{\partial s}\Big|_{s=0}\tau_{\pm}(s)\left(\left|\nabla_{g} f_{\pm}\right|^{2}_{g}+S_{g,u}\pm\frac{n}{2\tau_{\pm}}\right)
-\tau_{\pm}\bigg(h^{ij}\nabla_{i}\nabla_{j}f
+h^{ij}S_{ij}\\
&& \ \pm\frac{1}{2\tau_{\pm}}{\rm tr}_{g}h+4v\left(\langle du,df\rangle-\Delta_{g}u\right)\bigg)\bigg]
\frac{e^{-f_{\pm}}}{(4\pi\tau_{\pm})^{n/2}}dV_{g}.
\end{eqnarray*}
To prove the theorem, it is sufficient to show that
\begin{equation*}
\int_{M}\left(\left|\nabla_{g}f_{\pm}\right|^{2}_{g}+
S_{g,u}\pm\frac{n}{2\tau_{\pm}}\right)\frac{e^{-f_{\pm}}}{(4\pi\tau_{\pm})^{n/2}}dV=0.
\end{equation*}
Since $M$ is compact, we have
\begin{equation*}
0=\int_{M}\Delta_{g}\left(e^{-f_{\pm}}\right)=\int_{M}
\left(-\Delta_{g}f_{\pm}+\left|\nabla_{g} f_{\pm}\right|^{2}_{g}\right)e^{-f_{\pm}}dV.
\end{equation*}
Hence
\begin{eqnarray*}
&&\int_{M}\left(\left|\nabla_{g} f_{\pm}\right|^{2}+S_{g,u}\pm\frac{n}{2\tau_{\pm}}\right)\frac{e^{-f_{\pm}}}{(4\pi\tau_{\pm})^{n/2}}dV\\
&=&\int_{M}\left(2\Delta_{g}f_{\pm}-\left|\nabla_{g} f\right|^{2}_{g}+S_{g,u}\pm\frac{n}{2\sigma_{\pm}}\right)
\frac{e^{-f_{\pm}}}{(4\pi\tau_{\pm})^{n/2}}dV.
\end{eqnarray*}
Then, Lemma \ref{l10.1} now indicates
\begin{eqnarray*}
&&\int_{M}\left(\left|\nabla_{g} f_{\pm}\right|^{2}+S_{g,u}\pm\frac{n}{2\tau_{\pm}}\right)\frac{e^{-f_{\pm}}}{(4\pi\tau_{\pm})^{n/2}}dV\\
&=&\int_{M}\left(\frac{\pm f_{\pm}\mp n+\nu_{\pm}(g,u)}{\tau_{\pm}}\pm\frac{n}{2}\right)\frac{e^{-f_{\pm}}}{(4\pi\tau_{\pm})^{n/2}}dV\\
&=&\int_{M}\frac{1}{\tau_{\pm}}\left(\pm f_{\pm}\mp\frac{n}{2}+\nu_{\pm}(g,u)\right)\frac{e^{-f_{\pm}}}{(4\pi\tau_{\pm})^{n/2}}dV\\
&=&\frac{1}{\tau_{\pm}}\left(
\pm\frac{n}{2}-\nu_{\pm}(g,u)\mp\frac{n}{2}+\nu_{\pm}(g,u)\right) \ \ = \ \ 0.
\end{eqnarray*}
The sign $+$ corresponds to the gradient expanding soliton while $-$ to the gradient shrinker soliton.
\end{proof}

\begin{corollary} \label{c10.3}Suppose that $(M,g)$ is a compact Riemannian
manifold and $u$ a smooth function on $M$. Let $h$ be any symmetric covariant $2$-tensor on $M$ and set $g(s):=g+sh$. Let $v$ be any smooth function on $M$ and $u(s):=u+sv$. If $\nu_{\pm}(g(s),u(s))=\mathcal{W}_{\pm}(g(s),u(s),f_{\pm}(s),\tau_{\pm}(s))$ for some smooth function $f_{\pm}(s)$ with $\int_{M}e^{-f_{\pm}(s)}dV/(4\pi\tau_{\pm}(s))^{n/2}=1$ and a constant
$\tau_{\pm}(s)>0$, and $(g,u)$ is a critical point of
$\nu_{\pm}(\cdot,\cdot)$, then
\begin{equation*}
\mathcal{S}_{g,u}=\mp\frac{1}{2\tau_{\pm}}g, \ \ \ f_{\pm}\equiv{\rm constant}.
\end{equation*}
Thus, if $\mathcal{W}_{\pm}(g,u,\cdot,\cdot)$ achieve their minimum and $(g,u)$ is
a critical point of $\nu_{\pm}(\cdot,\cdot)$, then $(M,g,u)$ satisfies the static Einstein vacuum equation.
\end{corollary}

\begin{proof} According to Lemma \ref{l10.1} and Theorem \ref{t10.2}, we have
\begin{eqnarray*}
&&\tau_{\pm}\left(-2\Delta_{g}f_{\pm}+\left|\nabla_{g} f_{\pm}\right|^{2}_{g}-S_{g,u}\right)\pm f_{\pm}\mp n \ \ = \ \ -\nu_{\pm}\\
&=&-\int_{M}\left[\tau_{\pm}\left(S_{g,u}+\left|\nabla_{g} f\right|^{2}_{g}\right)
\mp f_{\pm}\pm n\right]\frac{e^{-f_{\pm}}}{(4\pi\tau_{\pm})^{n/2}}dV_{g},
\end{eqnarray*}
and hence
\begin{eqnarray*}
2\Delta_{g}f_{\pm}-\left|\nabla_{g} f_{\pm}\right|^{2}_{g}+S_{g,u}&=&\int_{M}\left(S_{g,u}
+\left|\nabla_{g} f_{\pm}\right|^{2}_{g}\right)\frac{e^{-f_{\pm}}}{(4\pi\tau_{\pm})^{n/2}}dV_{g}\\
&=&\int_{M}\left(S_{g,u}+\Delta_{g}f_{\pm}\right)\frac{e^{-f_{\pm}}}{(4\pi\tau_{\pm})^{n/2}}dV_{g}\\
&=&\mp\frac{n}{2\tau_{\pm}} \ \ = \ \ S_{g,u}+\Delta_{g}f_{\pm}.
\end{eqnarray*}
From this we get $\Delta_{g}f_{\pm}=|\nabla_{g} f_{\pm}|^{2}_{g}$. After taking the integration on both sides, the functions $f_{\pm}$ must be constant that imply $\mathcal{S}_{g}\pm\frac{1}{2\tau_{\pm}}g=0$.
\end{proof}

\begin{remark} \label{r10.4} In the situation of Corollary \ref{c10.3}, by normalization, we my choose $f_{\pm}=\frac{n}{2}$.
\end{remark}




\begin{thebibliography}{10}

\bibitem{CHI} Cao, H., Hamilton, R., Ilmanen, T., \textit{Gaussian densities and stability for some Ricci solitons}, arXiv: math.DG/0404165.

\bibitem{CZ} Cao, Huai-Dong; Zhu, Meng, \textit{On second variation of Perelman's Ricci shrinker entropy}, Math. Ann. {\bf 353}(2012), no. 3, 747--763. MR 2923948

\bibitem{C1} Cao, Xiaodong, \textit{Eigenvalues of
$\left(-\Delta+\frac{R}{2}\right)$ on manifolds with nonnegative
curvature operator}, Math. Ann., {\bf 337}(2007), no. 2, 435--441. MR 2262792 (2007 g: 53071)

\bibitem{C2} Cao, Xiaodong, \textit{First eigenvalues of geometric operators under the Ricci flow}, Proc. Amer. Math. Soc., {\bf 136}(2008), no. 11, 4075--4078. MR 2425749 (2009 f: 53098)

\bibitem{FIN} Feidman, Michael; Ilmanen, Tom;  Ni, Lei, \textit{Entropy and reduced distance for Ricci expanders}, J. Geom. Anal., {\bf 15}(2005), no. 1, 49--62. MR 2132265 (2006 b: 53091)


\bibitem{HHKL} He, Chun-Lei; Hu, Sen; Kong, De-Xing; Liu, Kefeng, \textit{Generalized Ricci flow I: Local existence and uniqueness}, Topology and physics, 151--171, Nankai Tracts Math., {\bf 12}, World Sci. Publ., Hackensack, NJ, 2008. MR 2503395 (2010 k: 53098)

\bibitem{Li} Li, Jun-Fang, \textit{Eigenvalues and energy functionals with monotonicity formulae under Ricci flow}, Math. Ann., {\bf 338}(2007), no. 4, 927--946. MR 2317755 (2008 c: 53068)

\bibitem{LY} Li, Yi, \textit{Generalized Ricci flow I: higher derivatives estimates for compact manifolds}, Analysis $\&$ PDE., {\bf 5}(2012), no. 4, 747--775.

\bibitem{L} List, B., \textit{Evolution of an extended Ricci flow system}, PhD thesis, AEI Potsdam, 2005.

\bibitem{M1} M\"uller, Reto, \textit{Monotone volume formulas for geometric flows}, J. Reine Angew. Math., {\bf 643}(2010), 39--57. MR 2658189 (2011 k: 53086)

\bibitem{M2} M\"uller, Reto, \textit{Ricci flow coupled with harmonic map flow}, Ann. Sci. \'Ec. Norm. Sup\'er. (4) {\bf 45}(2012), no. 1, 101--142. MR 2961788

\bibitem{OSW} Oliynyk, T.; Sunneeta, V.; Woolgar, E., \textit{A gradient flow for worldsheet nonlinear sigma models}, Nuclear Phys. B {\bf 739}(2006), no. 3, 441--458. MR 2214659 (2006 m: 81185)

\bibitem{S1} Streets, Jeffrey, \textit{Ricci Yang-Mills flow}, PH.D., Thesis, Duck University, 2007.

\bibitem{S2} Streets, Jeffrey, \textit{Regularity and expanding entropy for connection Ricci flow}, J. Geom. Phys., {\bf 58}(2008), no. 7, 900--912. MR 2426247 (2009 f: 53105)

\bibitem{S3} Streets, Jeffrey, \textit{Singularity of renormalization group flows}, J. Geom. Phys., {\bf 59}(2009), no. 1, 8--16. MR 2479257 (2010 a: 53143)

\bibitem{S4} Streets, Jeffrey, \textit{Ricci Yang-Mills flow on surfaces}, Adv. Math., {\bf 223}(2010), no. 2, 454--475. MR 2565538 (2011 c: 53164)

\bibitem{Y} Young, A., \textit{Modified Ricci flow on a principal bundle}, Ph.D., thesis, University of Texas at Austin, 2008.

\bibitem{Z} Zhu, Meng, \textit{The second variation of the Ricci expanding entropy}, Pacific J. Math., {\bf 251}(2011), no. 2, 499--510. MR 2811045 (2012 g: 53143)



\end{thebibliography}
\end{document}